\documentclass[preprint,sort&compress,final,10pt]{elsarticle}

\usepackage{amsmath}
\usepackage{amssymb}
\usepackage{graphicx}
\usepackage{latexsym}
\usepackage{epstopdf}
\usepackage{natbib}
\usepackage{color}
\usepackage{hyperref}
\usepackage{amsmath,amscd}

\newtheorem{theorem}{Theorem}[section]

\newtheorem{definition}[theorem]{Definition}
\newtheorem{lemma}[theorem]{Lemma}

\newtheorem{remark}[theorem]{Remark}
\numberwithin{equation}{section}

\def\Proof{\noindent{\bf Proof.}~}
\def\qed{\hfill$\square$\smallskip}

\def\Im{\mathrm{Im}}

\journal{\empty}
\date{}

\begin{document}

\begin{frontmatter}

\title{Invariant curves of almost periodic twist mappings}

\author[au1]{Peng Huang}

\address[au1]{School of Mathematics Sciences, Beijing Normal University, Beijing 100875, P.R. China.}

\ead[au1]{hp@mail.bnu.edu.cn}

\author[au1]{Xiong Li\footnote{Corresponding author. Partially supported by the NSFC (11571041) and the Fundamental Research Funds for the Central Universities.}}

\ead[au1]{xli@bnu.edu.cn}

\author[au2]{Bin Liu\footnote{Partially supported by the NSFC (11231001).}}

\address[au2]{School of Mathematical Sciences, Peking University, Beijing 100871, P.R. China.}

\ead[au3]{bliu@pku.edu.cn}

\begin{abstract}
In this paper we are concerned with the  existence of invariant curves of planar twist mappings which are almost periodic in a spatial variable. As an application of this result to differential equations we will discuss the existence of almost periodic solutions and the boundedness of all solutions for superlinear Duffing's equation with an almost periodic external force.
\end{abstract}


\end{frontmatter}

\section{Introduction}
In this paper we are  concerned with the  existence of invariant curves of the following planar almost periodic  mapping
\begin{equation}\label{M}
\mathfrak{M}:\quad \begin{array}{ll}
\left\{\begin{array}{ll}
x_1=x+y+f(x,y),\\[0.2cm]
y_1=y+g(x,y),
 \end{array}\right.\  (x,y)\in \mathbb{R} \times [a,b],
\end{array}
\end{equation}
where the perturbations $f(x,y)$ and $g(x,y)$ are almost periodic in $x$ with the frequency $\omega=(\cdots,\omega_\lambda,\cdots)_{\lambda\in \mathbb{Z}}$ and admit a rapidly converging Fourier series expansion, respectively.

When the mapping $\mathfrak{M}$ in (\ref{M}) is periodic and possesses the intersection property,\ Moser \cite{Moser62} considered  the twist mapping
$$
\mathfrak{M}_{0}:\quad \begin{array}{ll}
\left\{\begin{array}{ll}
x_1=x+\alpha(y)+\varphi_{1}(x,y),\\[0.2cm]
y_1=y+\varphi_{2}(x,y),
 \end{array}\right.
\end{array}
$$
where the perturbations $\varphi_{1},\varphi_{2}$ are assumed to be small and of periodic $2\pi$ in $x.$\ He obtained the existence of invariant closed curves of  $\mathfrak{M}_{0}$ which is of class $\mathcal{C}^{333}$.\ About  $\mathfrak{M}_{0}$,\ an analytic version of the invariant curve theorem was presented in \cite{Siegel97},\ a version in class $\mathcal{C}^{5}$ in  R\"{u}ssmann \cite{Russmann70} and a optimal version in class $\mathcal{C}^{p}$ with $ p>3$  in Herman \cite{Herman83, Herman86}. Ortega \cite{Ortega99}  obtained a variant of the small twist theorem and   also studied the existence of invariant curves of mappings with average small twist in  \cite{Ortega01}.

When the perturbations $f(x,y), g(x,y)$ in (\ref{M}) are quasi-periodic in $x$, there are some results about the existence of invariant curves of the following planar quasi-periodic  mapping
\begin{equation}\label{a7}
\mathfrak{M}_{1}:\quad \begin{array}{ll}
\left\{\begin{array}{ll}
x_1=x+\alpha+y +f(x,y),\\[0.2cm]
x_1=y+g(x,y),
 \end{array}\right.\ \ \ \   (x,y)\in \mathbb{R} \times [a,b],
\end{array}
\end{equation}
where the functions $f(x,y)$ and $g(x,y)$ are quasi-periodic in $x$ with the frequency $\omega=(\omega_1,\omega_2$,$\cdots,\omega_n)$, and $\alpha$ is a constant.

When the map $\mathfrak{M}_{1}$ in (\ref{a7}) is an exact symplectic map,\ $\omega_1,\omega_2,\cdots,\omega_n$, $2\pi\alpha^{-1}$ are sufficiently incommensurable and $f,g$ are real analytic in $x$ and $y$, Zharnitsky \cite{Zharnitsky00} proved the existence of invariant curves of the map $\mathfrak{M}_{1}$ and applied this result to present the boundedness of all solutions of Fermi-Ulam problem.\ His proof is based on the Lagrangian approach introduced by Moser \cite{Moser88} and used by Levi and Moser in \cite{Levi01} to show a proof of the twist theorem.

When the map $\mathfrak{M}_{1}$ in (\ref{a7}) is  reversible with respect to the involution $\mathcal{G}:(x,y)\mapsto (-x,y)$,\ that is,\ $\mathcal{G} \mathfrak{M}_1 \mathcal{G} = \mathfrak{M}_1^{-1}, \omega_1,\omega_2,\cdots,\omega_n,2\pi\alpha^{-1}$ satisfy the Diophantine condition
\begin{equation*}
\Big|{\langle k,\omega \rangle {\alpha \over {2\pi}}-j}\Big|\geq {\gamma \over {|k|^\tau}},\ \ \ \ \forall\ \  k \in \mathbb{Z}^n\backslash\{0\},\ \ \forall j \in \mathbb{Z},
\end{equation*}
and $f,g$ are real analytic in $x$ and $y$, Liu  \cite{Liu05} obtained some variants of the invariant curve theorem for quasi-periodic reversible mapping $\mathfrak{M}_{1}$. As an application,\ he used the invariant curve theorem  to investigate the existence of  quasi-periodic solutions and the boundedness of all solutions for an asymmetric oscillator depending quasi-periodically on time.

Recently,  in \cite{Huang1} we also considered  the  existence of  invariant curves of the planar quasi-periodic mapping $\mathfrak{M}$ in (\ref{M}). Instead of the exact symplecticity or reversibility assumption on $\mathfrak{M}$,\ this mapping $\mathfrak{M}$ is assumed to be possessing the intersection property,\ and  we obtained the  invariant curve theorem for  the quasi-periodic  mapping $\mathfrak{M}$ in the smooth case, other than analytic case, that is,  we assume that this mapping $\mathfrak{M}$ belongs to $\mathcal{C}^{p}$ with $p>2n+1$ and $n$ is the number of the frequency $\omega=(\omega_1,\omega_2,\cdots,\omega_n)$. We note that when $n=1$, quasi-periodic mappings are periodic mappings, and the optimal smoothness assumption is $\mathcal{C}^{p}$ with $p>3$. Hence our smoothness assumption for quasi-periodic mappings agrees with that for periodic mappings, and is optimal in this sense.

The application to an asymmetric oscillation about this theorem will be in a forthcoming paper \cite{Huang2}. However, the poincar\'{e} mapping of the asymmetric oscillation does not have the form as $\mathfrak{M}$, but has the following expression
\begin{equation}\label{f11}
\begin{array}{ll}
\mathcal{M}_{\delta}:\ \ \left\{\begin{array}{ll}
x_1=x+\beta+\delta l(x,y)+\delta f(x,y,\delta),\\[0.2cm]
y_1=y+\delta m(x,y)+\delta g(x,y,\delta),
\end{array}\right.\ \ \  (x,y)\in \mathbb{R} \times [a,b],
\end{array}
\end{equation}
where the functions $l,m, f, g$ are quasi-periodic in $x$ with the frequency $\omega=(\omega_1,\omega_2,\cdots,\omega_n),$
$f(x,y,0 ) =g(x,y, 0 ) = 0, $\ $\beta$ is a constant, $0<\delta< 1$ is a small parameter. Firstly, we \cite{Huang2} established the  existence of invariant curves of the  planar quasi-periodic  mapping $\mathcal{M}_{\delta}$ on the basis of the  invariant curve theorem obtained in \cite{Huang1} when the functions $l,m,f,g$ are smooth functions. After getting those invariant curves theorems, we study the existence of quasi-periodic solutions and the boundedness of all solutions for the asymmetric oscillation.

In this paper we focus on the almost periodic case, that is, the perturbations  $f(x,y), g(x,y)$ in (\ref{M}) are almost periodic in $x$ and admit a rapidly converging Fourier series expansion respectively,\ and we also  assume that the mapping $\mathfrak{M}$ in (\ref{M}) satisfies the intersection property,\ and  want to establish the invariant curve theorem for almost periodic mapping $\mathfrak{M}$ in the analytic case.\

After we get the invariant curve theorem, as an application, we shall study the existence of almost periodic solutions and the boundedness of all solutions for the following superlinear Duffing's equation
\begin{equation}\label{a1}
\ddot{x}+x^3=f(t),
\end{equation}
where $f(t)$ is a real analytic almost periodic function  with the frequency $\omega=(\cdots,\omega_\lambda,\cdots)$ and admits a rapidly converging Fourier series expansion.

It is well known that the longtime behaviour of a time dependent nonlinear differential equation
\begin{equation}\label{G}
\ddot{x}+f(t,x)=0,
\end{equation}
$f$ being periodic in $t$, can be very intricate. For example, there are equations having unbounded solutions but with infinitely many zeros and with nearby unbounded solution having randomly prescribed numbers of zeros and also periodic solution (see \cite{Dieckerhoff87}).

In contrast to such unboundedness phenomena one may look for conditions on the nonlinearity, in addition to the superlinear condition that
\begin{equation*}
{1\over x}f(t,x)\rightarrow \infty \ \ \ \text{as}\  |x|\rightarrow \infty,
\end{equation*}
which allow to conclude that all solutions of equation (\ref{G}) are bounded. For example, every solution of equation (\ref{a1})
with $p(t+1)=p(t)$ being continuous, is bounded. This result, prompted by Littlewood in \cite{Littlewood68}, is due to Morris \cite{Morris76}, who proved that there are infinitely many quasi-periodic solutions and the boundedness of all solutions of (\ref{a1}). In 1987, Dieckerhoff and Zehnder in \cite{Dieckerhoff87} extended this result to the general superlinear Duffing's equation. For recent development,\ we refer to\ \cite{{Liu98},  {Levi91}}\ and the references therein.

The rest of the paper is organized as follows.\ In Section 2,\ we first define real analytic almost periodic functions and their norms, then list some properties of them,\ and at last state the main invariant curve theorem  (Theorem \ref{thm2.11}) for the  almost periodic mapping $\mathfrak{M}$ given by (\ref{M}).\ The measure estimate  and the proof of Theorem \ref{thm2.11} are given in Sections 3, 4, 5 respectively. The small twist theorem is given in Section 6. In Section 7,\ we will prove the existence of almost periodic solutions and the boundedness of all solutions for superlinear Duffing's equation (\ref{a1}) with an almost periodic external force.

\section{Real analytic almost periodic functions and the main result}

\subsection{The space of  real analytic almost periodic functions}
Our aim is to find almost periodic solutions $x$ for superlinear Duffing's equation (\ref{a1}) which admit a rapidly converging Fourier series expansion, thus we have to  define  the space of a kind of real analytic almost periodic functions.

We first define  the space of real analytic quasi-periodic functions $Q(\omega)$ as in \cite[chapter 3]{Siegel97}, here the $n$-dimensional frequency vector $\omega=(\omega_1,\omega_2,\cdots,\omega_n)$ is rationally independent, that is, for any $k=(k_1,k_2,\cdots,k_n) \neq 0$,\ $\langle k,\omega \rangle =\sum k_j \omega_j \neq 0$.

\begin{definition}\label{def2.1}
A function $f:\mathbb{R} \rightarrow \mathbb{R}$ is called real analytic quasi-periodic with the frequency $\omega=(\omega_1,\omega_2,\cdots,\omega_n)$,  if there exists a real analytic function
$$F: \theta=(\theta_1,\theta_2,\cdots,\theta_n) \in \mathbb{R}^n \rightarrow \mathbb{R}$$
such that $f(t)=F(\omega_1t, \omega_2t,\cdots, \omega_nt)\ \text{for all}\ t\in \mathbb{R}$, where $F$ is $2\pi$-periodic in each variable and bounded in a complex neighborhood $\Pi_{r}^{n}=\{(\theta_{1},\theta_{2},\cdots,\theta_{n})\in \mathbb{C}^n : |\Im\ \theta_{j}|\leq r, j=1,2,\cdots, n \}$  of\, $\mathbb{R}^n$ for some $r > 0$. Here we call $F(\theta)$ the shell function of $f(t)$.
\end{definition}

Denote by $Q(\omega)$  the set of real analytic quasi-periodic functions with the frequency $\omega=(\omega_1,\omega_2,\cdots,\omega_n)$.  Given $f(t)\in Q(\omega)$, the shell function $F(\theta)$ of  $f(t)$ admits a Fourier series expansion
$$F(\theta)=\sum \limits_{k\in \mathbb{Z}^n} f_{k}e^{i \langle k,\theta \rangle },$$
where $k=(k_1,k_2,\cdots,k_n)$,  $k_j$ range over all integers and the coefficients $f_{k}$ decay exponentially with $|k|=|k_1|+|k_2|+\cdots+|k_n|$, then $f(t)$ can be represented as a Fourier series of the type from the definition,
$$
\begin{array}{ll}
f(t)=\sum \limits_{k\in \mathbb{Z}^n} f_{k}e^{i \langle k,\omega \rangle t}.
\end{array}
$$

In the following we define the norm of the real analytic quasi-periodic function $f(t)$ through that of the corresponding shell function $F$.
\begin{definition}
For $r>0$, let $Q_{r}(\omega)\subseteq Q(\omega)$ be the set of real analytic quasi-periodic functions $f$ such that the corresponding shell functions $F$ are bounded on the subset $\Pi_{r}^{n}$\ with the supremum norm
$$\big|F\big|_{r}=\sup \limits_{\theta\in \Pi_{r}^{n}}|F(\theta)|=\sup \limits_{\theta\in \Pi_{r}^{n}}\Big|\sum_{k}f_{k}e^{i\langle k,\theta\rangle}\Big|<+\infty.$$
Thus we define $\big|f\big|_{r}:=\big|F\big|_{r}.$
\end{definition}

Similarly, one can give the definition of real analytic almost periodic functions with the frequency $\omega=(\cdots,\omega_\lambda,\cdots)_{\lambda\in \mathbb{Z}}$, which is a bilateral infinite sequence of rationally independent frequency, that is to say, any finite segments of $\omega$  are rationally independent. For this purpose, we first define analytic functions on some infinite dimensional space (see \cite{Dineen99}).
\begin{definition}
Let $X$ be a complex Banach space. A function $f : U\subseteq X \rightarrow \mathbb{C}$, where $U$ is an open subset of $X$, is called analytic if $f$ is continuous on $U$, and $f|_{U\cap X_1}$ is analytic in the classical sense as a function of several complex variables for each finite dimensional subspace $X_1$ of $X$.
\end{definition}

Note that for the bilateral infinite sequence of rationally independent frequency $\omega=(\cdots,\omega_\lambda,\cdots)_{\lambda\in \mathbb{Z}}$,
$$\langle k,\omega \rangle=\sum \limits_{\lambda\in \mathbb{Z}} k_\lambda \omega_\lambda,$$
where due to the spatial structure of the perturbation $k$ runs over integer vectors whose support
$$\mbox{supp}\, k=\big\{\lambda\ : \ k_{\lambda}\neq 0\big\}$$
is a finite set of $\mathbb{Z}$.

\begin{definition}\label{def2.2}
A function $f:\mathbb{R} \rightarrow \mathbb{R}$ is called real analytic almost periodic with the frequency $\omega=(\cdots,\omega_\lambda,\cdots)\in \mathbb{R}^\mathbb{Z}$, if there exists a real analytic function
$$F: \theta=(\cdots,\theta_\lambda,\cdots) \in \mathbb{R}^{\mathbb{Z}} \rightarrow \mathbb{R},$$
which admit a rapidly converging Fourier series expansion
$$F(\theta)=\sum \limits _{A\in\mathcal{S}} F_{A}(\theta),$$
where
$$F_{A}(\theta)={\sum \limits _{\mbox{supp}\, k \subseteq A}} f_{k}\,e^{i \langle k,\theta \rangle },$$
and $\mathcal{S}$ is a family of  finite subsets $A$ of $\mathbb{Z}$ with $\mathbb{Z}\subseteq \bigcup\limits_{A\in\mathcal{S}} A, \langle k,\theta \rangle=\sum \limits_{\lambda\in \mathbb{Z}} k_\lambda \theta_\lambda,$
such that $f(t)=F(\omega t)\ \text{for all}\ t\in \mathbb{R}$, where $F$ is $2\pi$-periodic in each variable and bounded in a complex neighborhood  $\Pi_{r}=\Big\{\theta=(\cdots,\theta_\lambda,\cdots) \in \mathbb{C}^{\mathbb{Z}} :\ |\Im\,\theta|_{\infty}\leq r\Big\}$ for some $r > 0$, where $|\Im\,\theta|_{\infty}=\sup \limits_{\lambda\in \mathbb{Z}} |\Im\,\theta_{\lambda}|$. Here $F(\theta)$ is called the shell function of $f(t)$.
\end{definition}

Suppose that the function $f(t)$ has the Fourier exponents $\{\Lambda_\lambda: \lambda\in \mathbb{Z}\}$, and its basis is $\{\omega_\lambda: \lambda\in \mathbb{Z}\}$. Then for any $\lambda\in \mathbb{Z}$, $\Lambda_\lambda$ can be expressed into
$$
\Lambda_\lambda=r_{\lambda_1}\omega_{\lambda_1}+\cdots+r_{\lambda_{j(\lambda)}}\omega_{\lambda_{j(\lambda)}},
$$
where $r_{\lambda_1},\cdots, r_{\lambda_{j(\lambda)}}$ are rational numbers. Therefore,
$$
\mathcal{S}=\{(\lambda_1,\cdots,\lambda_{j(\lambda)}):\lambda\in\mathbb{Z}\}.
$$
Thus, this family  $\mathcal{S}$ is not totally arbitrary.\ Rather,\ $\mathcal{S}$ has to be a spatial structure on $\mathbb{Z}$ characterized by the property that the union of any two sets in $\mathcal{S}$ is again in $\mathcal{S}$,\ if they intersect :
$$A,\ B \in \mathcal{S},\ \ \ A\cap B \neq \varnothing\ \ \ \Rightarrow \ \ \ A \cup B \in \mathcal{S}.$$
Moreover, we define
\begin{equation}\label{a00000}
{\mathbb{Z}}_{0}^{\mathbb{Z}}:=\Big\{k=(\cdots,k_\lambda,\cdots)\in\mathbb{Z}^{\mathbb{Z}} :\ \mbox{supp}\,k\subseteq A,\ A\in\mathcal{S}\Big\}.
\end{equation}

Denote by $AP(\omega)$  the set of real analytic almost periodic functions with the frequency $\omega$ which are given by Definition \ref{def2.2}. From the definition, given $f(t)\in AP(\omega)$, the shell function $F(\theta)$ admits a rapidly converging Fourier series expansion
$$F(\theta)=\sum \limits _{A\in\mathcal{S}} F_{A}(\theta),$$
where
$$F_{A}(\theta)={\sum \limits _{\mbox{supp}\, k \subseteq A}} f_{k}\,e^{i \langle k,\theta \rangle }.$$
From the definitions of the support $\mbox{supp}\, k$ and $A$, we know that $F_{A}(\theta)$ is real analytic and $2\pi$-periodic in $\big\{\theta_\lambda\ :\ \lambda\in A\big\}$.

As a consequence of the definitions of $\mathcal{S}$ and the support $\mbox{supp}\, k$ of $k$ and $\mathbb{Z}_{0}^{\mathbb{Z}}$,\ the Fourier series of the shell function $F(\theta)$ has another form
$$F(\theta)=\sum \limits _{k \in {{\mathbb{Z}}_{0}^{\mathbb{Z}} }} f_{k}e^{i \langle k,\theta \rangle },$$
then $f(t)$ can be represented as a Fourier series of the type
\begin{equation*}
f(t)={\sum \limits_{A\in \mathcal{S}}}\ {\sum \limits _{\mbox{supp}\, k \subseteq A}} f_{k}e^{i \langle k,\omega \rangle t},\quad\quad \langle k,\omega \rangle=\sum \limits_{\lambda\in \mathbb{Z}} k_\lambda \omega_\lambda
\end{equation*}
or
\begin{equation*}
f(t)=\sum \limits _{k \in {{\mathbb{Z}}_{0}^{\mathbb{Z}} }}f_{k}e^{i \langle k,\omega \rangle t}.
\end{equation*}
If we define
$$
f_{A}(t)={\sum \limits _{\mbox{supp}\, k \subseteq A}} f_{k}\,e^{i \langle k,\omega \rangle t},
$$
then
$$f(t)=\sum \limits _{A\in\mathcal{S}} f_{A}(t).$$
From the definitions of the support $\mbox{supp}\, k$ and $A$, we know that $f_{A}(t)$ is a real analytic quasi-periodic function with the frequency $\omega_A=\big\{\omega_\lambda\ :\ \lambda\in A\big\}$.

\subsection{The norms of real analytic almost periodic functions}

In the following we will give a kind of the norm for real analytic almost periodic functions. Before we describe the norm,  some more definitions  and  notations  are useful.

The main ingredient of our perturbation theory is a nonnegative weight function
$$[\,\cdot\,]\ :\ \ \ A\ \mapsto\ [A]$$
defined on $\mathcal{S}.$ The weight of a subset may reflect its
size,\ its location or something else.\ This, however,\ is immaterial for the perturbation theory itself.\ Here only the properties of monotonicity and subadditivity are required :
$$A\subseteq B \ \ \ \ \Rightarrow\ \ \ \  [A]\leq [B],$$
$$A\cap B \neq \varnothing \ \ \ \ \Rightarrow \ \ \ \ [A\cup B]+[A\cap B]\leq [A]+[B]$$
for all $A,\ B$ in $\mathcal{S}$.\ Throughout this paper, we always use the following weight function
$$[A]=1+\sum\limits_{i\in A}log^{\varrho}(1+|i|),$$
where $\varrho>2$ is a constant.

Now we can define the norm of the real analytic almost periodic function $f(t)$ through that of the corresponding shell function $F$ just like in the quasi-periodic case.

\begin{definition}
Let $AP_{r}(\omega)\subseteq AP(\omega)$ be the set of real analytic almost periodic functions $f$ such that the corresponding shell functions $F$ are bounded on the subset $\Pi_{r}$\ with the norm
$$\|F\|_{m,r}=\sum \limits _{A\in\mathcal{S}} |F_{A}|_{r}\, e^{m[A]}=\sum \limits _{A\in\mathcal{S}} |f_{A}|_{r}\, e^{m[A]}<+\infty,$$
where $m>0$ is a constant and
$$|F_{A}|_{r}=\sup\limits_{\theta\in\Pi_{r}}|F_{A}(\theta)|=\sup\limits_{\theta\in\Pi_{r}} \Bigg|\sum \limits _{\mbox{supp}\, k \subseteq A}f_{k}\, e^{i\langle k,\theta\rangle}\Bigg|=|f_{A}|_{r}.$$
Hence we define
$$\|f\|_{m,r}:=\|F\|_{m,r}.$$

If $f(\cdot,y)\in AP_r(\omega)$, and the corresponding shell functions $F(\theta,y)$ are real
analytic in the domain $D(r,s)=\big\{(\theta,y)\in \mathbb{C}^\mathbb{Z} \times \mathbb{C}\ :\ |\Im\,\theta|_{\infty}<r, |y-\alpha|<s\big\}$ with some $\alpha\in \mathbb{R}$, we define
$$\|f\|_{m,r,s}=\sum \limits _{A\in\mathcal{S}} |F_{A}|_{r,s}\, e^{m[A]}=\sum \limits _{A\in\mathcal{S}} |f_{A}|_{r,s}\, e^{m[A]},$$
where
$$|F_{A}|_{r,s}=\sup\limits_{(\theta,y)\in D(r,s)}|F_{A}(\theta,y)|=\sup\limits_{(\theta,y)\in D(r,s)} \Bigg|\sum \limits _{\mbox{supp}\, k \subseteq A}f_{k}(y)\, e^{i\langle k,\theta\rangle}\Bigg|=|f_{A}|_{r,s}.$$
\end{definition}

\subsection{The properties of  real analytic almost periodic functions}
In the following some properties of real analytic almost periodic functions are given.
\begin{lemma}\label{lem2.4}
The following statements are true:\\
$(i)$ Let $f(t),g(t)\in AP(\omega)$,\ then $f(t)\pm g(t), g(t+f(t))\in AP(\omega);$\\[0.1cm]
$(ii)$ Let $f(t)\in AP(\omega)$ and $\tau=\beta t +f(t)\ (\beta+f'>0)$,\ then the inverse relation is given by $t={\beta^{-1}}\tau +g(\tau)$  and $g\in AP({\omega / \beta})$.\ In particular,\ if $\beta=1$,\ then $g\in AP(\omega).$
\end{lemma}

\Proof One can use a similar way in \cite[chapter 3]{Siegel97} to prove this lemma.\\
$(i)$ It is easy to see that $f(t)\pm g(t)\in AP(\omega)$. Now we prove $g(t+f(t))\in AP(\omega)$.  We claim that if the shell functions of $f(t),g(t)$ are $F(\theta),G(\theta)$ respectively, from the definitions of the shell function and $\omega$, we know that $\omega F(\theta)$ and $\theta$ is of the same dimension, then the shell function of $g(t+f(t))$ is $G(\theta+\omega F(\theta)).$ Indeed, since $f(t),g(t)\in AP(\omega)$, then we have
$$f(t)=\sum \limits _{k \in {{\mathbb{Z}}_{0}^{\mathbb{Z}} }}f_{k}e^{i \langle k,\omega \rangle t},\  F(\theta)=\sum \limits _{k \in {{\mathbb{Z}}_{0}^{\mathbb{Z}} }} f_{k}e^{i \langle k,\theta \rangle }$$
and
$$g(t)=\sum \limits _{k \in {{\mathbb{Z}}_{0}^{\mathbb{Z}} }} g_{k}e^{i \langle k,\omega \rangle t},\ G(\theta)=\sum \limits _{k \in {{\mathbb{Z}}_{0}^{\mathbb{Z}} }} g_{k}e^{i \langle k,\theta \rangle }.$$
Hence,
\begin{eqnarray*}
g(t+f(t))&=&\sum \limits _{k \in {{\mathbb{Z}}_{0}^{\mathbb{Z}} }}g_{k}e^{i \langle k,\omega \rangle (t+f(t))}
=\sum \limits _{k \in {{\mathbb{Z}}_{0}^{\mathbb{Z}} }} g_{k}e^{i \langle k,\omega \rangle t+ i \langle k,\omega \rangle f(t)}\\[0.0001cm]
&=&\sum \limits _{k \in {{\mathbb{Z}}_{0}^{\mathbb{Z}} }} g_{k} e^{{i \langle k,\omega \rangle t+ i \langle k,\omega \rangle \sum \limits _{k \in {{\mathbb{Z}}_{0}^{\mathbb{Z}} }} f_{k}e^{i \langle k,\omega \rangle t}}}.
\end{eqnarray*}
From the definition of the shell function, we know that
\begin{align*}
\sum \limits _{k \in {{\mathbb{Z}}_{0}^{\mathbb{Z}} }} g_{k}e^{{i \langle k,\theta \rangle + i \langle k,\omega \rangle \sum \limits _{k \in {{\mathbb{Z}}_{0}^{\mathbb{Z}} }} f_{k}e^{i \langle k,\theta\rangle }}}
&=\sum \limits _{k \in {{\mathbb{Z}}_{0}^{\mathbb{Z}} }} g_{k}e^{i \langle k,\theta \rangle + i \langle k,\omega \rangle F(\theta)}\\[0.0001cm]
&=\sum \limits _{k \in {{\mathbb{Z}}_{0}^{\mathbb{Z}} }} g_{k}e^{i \langle k,\theta+\omega F(\theta) \rangle }=G(\theta+\omega F(\theta))
\end{align*}
is the shell function of $g(t+f(t))$.\\
$(ii)$ To prove this assertion, we may assume that $\beta=1$,\ or else replace $\tau$ by $\beta \tau$.
If $f(t)$ is represented by its shell function $F(\theta)$ and the unknown function $g(\tau)$ by $G(\theta)$, where $g(\tau)=G(\omega \tau)$, then
\begin{equation}\label{b4}
G(\theta)+F(\theta+\omega G)=0.
\end{equation}
We replace this equation by
\begin{equation*}
G+\sigma F(\theta+\omega G)=0
\end{equation*}
and seek a solution $G=G(\theta;\sigma)$ for $0\leq\sigma\leq 1$ with period $2\pi$ in each of the variables $\theta_{\lambda}\,( \lambda\in \mathbb{Z} )$. Differentiating this last equation  with respect to $\sigma$, we are led to the differential equation
\begin{equation}\label{b5}
{{\partial G}\over{\partial\sigma}}=\phi(\theta+\omega G,\sigma),\ \ G(\theta,0)=0,
\end{equation}
where
$$\phi(\theta,\sigma)=-F(\theta)\bigg(1+\sigma\sum \limits_{\lambda\in \mathbb{Z}}\omega_\lambda F_{\theta_\lambda}\bigg)^{-1}.$$
By the assumption on $\beta+f'$, the denominator on the right of the last expression is bounded away from $0$ when $\sigma=1$ and $\theta_\lambda=\omega_\lambda t \,(\lambda\in\mathbb{Z})$, and is indeed positive. On the other hand, the
vectors with components $\omega_\lambda t+2 k_\lambda \pi\,(\lambda\in\mathbb{Z})$ for integer $k_\lambda$ and real $t$ are dense in $\mathbb{Z}$-dimensional Euclidean space, and therefore the denominator is actually positive and bounded away from $0$ for all real $\theta$ and $0\leq\sigma\leq 1$. It follows that $\phi(\theta,\sigma)$ is real analytic and of period $2\pi$ in $\theta_{\lambda}$.

The solution $G(\theta,\sigma)$ of (\ref{b5}) is now constructed by means of the
standard existence theorem for ordinary differential equations. To show that $G(\theta,\sigma)$ is analytic for $|\Im\,\theta|_{\infty}$ sufficiently small and for $0\leq\sigma\leq 1$, it suffices to verify that as we continue the solution we do not leave the region of analyticity of the differential equation. Assuming that $\phi$ is analytic in $|\Im\,\theta|_{\infty}<\tilde{\delta}$, $0\leq\sigma\leq 1$, where it satisfies $\sum\limits_{\lambda\in\mathbb{Z}}|\phi_{\theta_\lambda}|<M,$ we set
$$\tilde{\varrho}=\tilde{\delta} e^{-M|\omega|_{\infty}},$$
where $|\omega|_{\infty}=\sup\limits_{\lambda\in \mathbb{Z}}|\omega_\lambda|<+\infty$, and claim that in the region $|\Im\,\theta|_{\infty}<\tilde{\varrho}$, $0\leq\sigma\leq 1$, the function $G(\theta,\sigma)$ satisfies the estimate $\tilde{\varrho}+|\omega|_{\infty}|\Im\,G|<\tilde{\delta}$. This, of course, will imply that the solution exists and is analytic in this region. We now verify the claim. Suppose the
estimate $\tilde{\varrho}+|\omega|_{\infty}|\Im\,G|<\tilde{\delta}$ does not hold over the entire interval $0\leq\sigma\leq 1$ as we continue the solution. Since it certainly
holds initially, there would exist a smallest number  $0<\sigma^*\leq 1$ for which
$$\sup\limits_{\theta} |\Im\,G|\geq{{\tilde{\delta}-\tilde{\varrho}}\over{|\omega|_{\infty}}}\ \ \ (\sigma=\sigma^*),$$
which means that the solution $G(\theta,\sigma)$ of (\ref{b5}) exists for $0\leq\sigma<\sigma^*.$

Since $\phi$ is real analytic, we have $\overline{\phi(\theta,\sigma)}=\phi(\bar{\theta},\sigma)$, and therefore
\begin{eqnarray*}
\Bigg|{\partial\over{\partial\sigma}}\Im\,G\Bigg|&=&{1\over 2}\Bigg|{\partial\over{\partial\sigma}}(G-\bar{G})\Bigg|\\[0.2cm]
&=& {1\over 2}\Big|\phi(\theta+\omega G,\sigma)-\phi(\bar{\theta}+\omega \bar{G},\sigma)\Big|\\[0.2cm]
&\leq& M\sup\limits_{\lambda\in\Lambda}|\Im\,\theta_\lambda+\omega_\lambda\Im\,G|\\[0.2cm]
&<&  M(\tilde{\varrho}+|\omega|_{\infty}|\Im\,G|)
\end{eqnarray*}
is valid for $0\leq\sigma<\sigma^*$. By
a comparison argument we would then conclude that
$$|\Im\,G|<h(\sigma)\ \ \ (0<\sigma\leq\sigma^*),$$
where $h(\sigma)$ is the solution of the equation
$${{dh}\over{d\sigma}}=M(\tilde{\varrho}+|\omega|_{\infty}h), \ \ \ h(0)=0.$$
And since
$$h(\sigma^*)={{\tilde{\varrho}}\over{|\omega|_{\infty}}}\Big(e^{M|\omega|_{\infty} \sigma^*}-1\Big)\leq{{\tilde{\varrho}}\over{|\omega|_{\infty}}}\Big(e^{M|\omega|_{\infty} }-1\Big)={{\tilde{\delta}-\tilde{\varrho}}\over{|\omega|_{\infty}}},$$
we would have
$$|\Im\,G|<{{\tilde{\delta}-\tilde{\varrho}}\over{|\omega|_{\infty}}}\ \ \ (\sigma=\sigma^*),$$
contrary to the choice of $\sigma^*$. This verifies our claim, and we can conclude that the solution $G(\theta,\sigma)$ of (\ref{b5}) exists for $0\leq\sigma\leq 1$, is real analytic, and by the uniqueness theorem, has period $2\pi$ in the $\theta_{\lambda}$. It is easy to verify that $G=G(\theta,1)$ is the solution of our original equation (\ref{b4}). Thus $g(\tau)=G(\omega\tau,1)$ belongs to $AP(\omega)$, and the assertion about the form of the inverse function is proved. \qed

\subsection{The frequency of the perturbations}
In this paper, the frequency $\omega=(\cdots,\omega_\lambda,\cdots)$ of the perturbations is not only rationally independent with $|\omega|_{\infty}=\sup\{|\omega_\lambda|:\lambda\in \mathbb{Z}\}<+\infty$,  but also it satisfies the strongly nonresonant  condition.

In a crucial fashion the weight function determines the nonresonance conditions for the small divisors arising in this theory. As we will do later on by way of an appropriate norm,\ it   suffices to estimate these small divisors from below not only in terms of the norm of $k$
$$|k|=\sum \limits_{\lambda \in \mathbb{Z}} |k_{\lambda}|,$$
but also in terms of the weight of its support
$$[[k]]=\min\limits_{\mbox{supp}\, k \subseteq A \in \mathcal{S}}[A].$$
Then the nonresonance conditions read
\begin{equation}\label{b3}
\begin{array}{ll}
{{| \langle k,\omega \rangle|} \geq {c \over {\Delta([[k]])\Delta(|k|)}}},\ \ \ \ 0\neq k \in \mathbb{Z}_{0}^{{\mathbb{Z}}},
\end{array}
\end{equation}
where, as usual, $c$ is a positive parameter and $\Delta$ some fixed approximation function as described in the following. One and the same approximation function is taken here in both places for simplicity, since the generalization is straightforward. A nondecreasing function $\Delta\ :\ [1,\infty) \rightarrow [1,\infty)$ is called an approximation function,\ if
\begin{equation}\label{b15}
 {{\log\Delta(t)}\over t} \searrow 0,\ \ \ \ 1\leq t \rightarrow\infty,
\end{equation}
and
$$\int_{1}^{\infty} {{\log\Delta(t)}\over {t^2}}\, dt <\infty.$$
In addition,\ the normalization $\Delta(1)=1$ is imposed for definiteness.

In the following we will give a criterion for the existence of strongly nonresonant
frequencies. It is based on growth conditions on the distribution function
$$N_n(t)=\mbox{card}\ \big\{A\in\mathcal{S}\ :\ \mbox{card}(A)=n,\ [A]\leq t\big\}$$
for $n\geq 1$ and $t\geq 0$.

\begin{lemma}\label{lem2.11}
There exist a constant $N_0$ and an approximation function $\Phi$ such that
$$
\begin{array}{ll}
N_n(t)\leq\left\{\begin{array}{ll}
0,\ \ \ \ \ \ \ \ \ \  t<t_n\\[0.4cm]
N_0 \Phi(t),\ \ \ t\geq t_n
 \end{array}\right.
\end{array}
$$
with a sequence of real numbers $t_n$ satisfying
$$n\log^{\varrho-1} n\leq t_n\sim n\, \log^{\varrho} n$$
for $n$ large with some exponent $\varrho-1>1$, we say $a_n \sim b_n$, if there are two constants $c, C$ such that $c a_n\le b_n\le C a_n$ and $c, C$ are independent of $n$.
\end{lemma}
\Proof We note that the weight function is
$$[A]=1+\sum\limits_{i\in A}\log^{\varrho}(1+|i|),\ \ \ \varrho>2.$$
Let $B_m=\{i : |i|\leq m\},$ then we have $[B_m] \sim |B_m|\log^{\varrho}|B_m|$, where $|B_m| = \text{card}(B_m)$, and $B_m$ has the lowest weight among all sets with the same
number of elements. It follows that
$$N_n(t)=0\ \ \  \text{for}\ \ \ t\leq t_n \sim n\log^{\varrho}n.$$

Next, let $t\geq 0$ be arbitrary, and consider the collection of all sets $A$ with $n$ elements and weight not bigger than $t$. Picking any element from $A$ with weight $0\leq t-s\leq t$, the remaining $n-1$ elements have total weight not bigger than $s$. This leads to the estimate
$$N_n(t)\leq{1\over n}\int_{0}^{t}W(t-s)dN_{n-1}(s),$$
where $W$ is any continuous function bounding $N_1$ from above. Integrating by parts
and assuming that $W(0)=0$ the role of $W$ and $N_{n-1}$ can be interchanged. And
proceeding by induction, we obtain
$$N_n(t)\leq{1\over {n!}}\int\limits_{t_1+\cdots+t_n\leq t}dW(t_1)\cdots dW(t_n).$$
Now choose $W(t)=\kappa e^{t^\mu}-\kappa$ with $\mu=\varrho^{-1}$ and a
suitable constant $\kappa\geq 2$. Then
$$N_n(t)\leq{{\kappa^n}\over {n!}}\int\limits_{t_1+\cdots+t_n\leq t}\exp(t_1^\mu+\cdots+t_n^\mu)dt_1^\mu\ldots dt_n^\mu.$$
On the domain of integration, $t_1^\mu+\cdots+t_n^\mu\leq n^{1-\mu}(t_1+\cdots+t_n)^\mu\leq n^{1-\mu}t^\mu$, while
the integral of $dt_1^\mu\ldots dt_n^\mu$ over $[0,t]^n$ is bounded by $t^{\mu n}$. Hence,
$$N_n(t)\leq{{\kappa^n t^{\mu n}}\over {n!}}\exp(n^{1-\mu}t^\mu)\leq \exp(\kappa t^\mu)\exp(n^{1-\mu}t^\mu).$$
Finally, to eliminate the dependence on $n$ for large $n$, recall that $t\geq t_n \sim n\log^{\varrho}n.$ Hence, $n \sim t_n\log^{\varrho}t_n\leq t/\log^{\varrho}t$ for $t\geq t_n$ and so
$$n^{1-\mu}t^\mu\lesssim {t\over{\log^{\varrho-1}t}},\ \ \ \varrho(1-\mu)=\varrho-1>1.$$
This shows that $N_n(t)$ is bounded from above by a constant multiple of a fixed
approximation function independently of $n$ as required by Lemma \ref{lem2.11}.\qed

According to Lemma \ref{lem2.11}, there exist an approximation function $\Delta_0$ and a probability measure $\mu$ on the parameter space $\mathbb{R}^{\mathbb{Z}}$ with support at any prescribed point such that the  measure of the set of $\omega$ satisfying the following inequalities
\begin{equation}\label{b6}
\begin{array}{ll}
{{| \langle k,\omega \rangle|} \geq {c \over {\Delta_0([[k]])\Delta_0(|k|)}}},\ \ \ \ c>0,\ \forall\ 0\neq k \in \mathbb{Z}_{0}^{{\mathbb{Z}}}
\end{array}
\end{equation}
is positive for a suitably small $c$,\ the proof can be found in \cite{Poschel90},\ we omit it here.

Throughout this paper,\ we always assume that the frequency  $\omega=(\cdots,\omega_\lambda,\cdots)$  satisfies the nonresonance condition (\ref{b6}).

\subsection{The main result}
We choose a rotation number $\alpha$ satisfying the inequalities
\begin{equation}\label{b7}
\begin{array}{ll}
\left\{\begin{array}{ll}
a+\gamma\leq \alpha\leq b-\gamma,\\[0.4cm]
\Big|{\langle k,\omega \rangle {\alpha \over {2\pi}}-j}\Big|\geq {\gamma \over {\Delta([[k]])\Delta(|k|)}},\ \ \ \  \mbox{for all}\ \  k \in \mathbb{Z}_{0}^{{\mathbb{Z}}}\backslash\{0\},\ \  j \in \mathbb{Z}
 \end{array}\right.
\end{array}
\end{equation}
with some positive constant $\gamma$, where $\Delta$ is some approximation function (see Theorem \ref{thm2.12} below).

\begin{definition}\label{def2.4}
Let $\mathfrak{M}$ be the mapping given by (\ref{M}). It is said that $\mathfrak{M}$ has the intersection property if
$$\mathfrak{M}(\mathbf{\Gamma}) \cap \mathbf{\Gamma} \neq \emptyset$$
for every  curve $\mathbf{\Gamma}:x=\xi+\varphi(\xi),\ y=\psi(\xi)$,\ where $\varphi$ and $\psi$ are real analytic and almost periodic in $\xi$ with the frequency  $\omega=(\cdots,\omega_\lambda,\cdots)$.
\end{definition}

\begin{definition}\label{def2.400}
Let $\mathfrak{M}$ be the mapping given by (\ref{M}). We say that $\mathfrak{M}:\mathbb{R} \times [a,b]\to \mathbb{R}^2$  is exact symplectic if
$\mathfrak{M}$ is  symplectic with respect to the usual symplectic structure $dy\wedge dx$ and for every  curve $\mathbf{\Gamma}:x=\xi+\varphi(\xi),\ y=\psi(\xi)$,\ where $\varphi$ and $\psi$ are real analytic and almost periodic in $\xi$ with the frequency  $\omega=(\cdots,\omega_\lambda,\cdots)$, we have
$$
\lim_{T\to+\infty}\frac{1}{2T}\int_{-T}^T\, ydx=\lim_{T\to+\infty}\frac{1}{2T}\int_{-T}^T\, y_1dx_1.
$$
\end{definition}

We claim that if the mapping $\mathfrak{M}$ is an exact symplectic map, then it has intersection property. In order to  prove this result, we first give an useful  lemma and some properties of generating functions of almost periodic monotonic twist maps.

\begin{lemma}\label{lem 7.1}
If $y = y(x)$ is real analytic almost periodic   in $x$  with the frequency  $\omega=(\cdots,\omega_\lambda,\cdots)$ and $F(x,y)$ is real analytic almost periodic  in $x$ with the same frequency, then $F(x,y(x))$ is also real analytic almost periodic in $x$ with the same frequency.
\end{lemma}
\Proof We claim that the shell functions of  $y(x),F(x,y)$ are $Y(\theta),\widetilde{F}(\theta,y)$ respectively, then the shell function of $F(x,y(x))$ is $\widetilde{F}(\theta,Y(\theta)).$ Indeed, since $y(x),F(x,y)\in AP(\omega)$, then we have
$$y(x)=\sum \limits _{k \in {{\mathbb{Z}}_{0}^{\mathbb{Z}} }}y_{k}e^{i \langle k,\omega \rangle x},\  Y(\theta)=\sum \limits _{k \in {{\mathbb{Z}}_{0}^{\mathbb{Z}} }} y_{k}e^{i \langle k,\theta \rangle }$$
and
$$F(x,y)=\sum \limits _{k \in {{\mathbb{Z}}_{0}^{\mathbb{Z}} }} F_{k}(y)e^{i \langle k,\omega \rangle x},\ \widetilde{F}(\theta,y)=\sum \limits _{k \in {{\mathbb{Z}}_{0}^{\mathbb{Z}} }} F_{k}(y)e^{i \langle k,\theta \rangle }.$$
Hence,
\begin{eqnarray*}
F(x,y(x))=\sum \limits _{k \in {{\mathbb{Z}}_{0}^{\mathbb{Z}} }} F_{k}\Bigg(\sum \limits _{k \in {{\mathbb{Z}}_{0}^{\mathbb{Z}} }}y_{k}e^{i \langle k,\omega \rangle x}\Bigg)e^{i \langle k,\omega \rangle x}.
\end{eqnarray*}
From the definition of the shell function, we know that
\begin{align*}
\sum \limits _{k \in {{\mathbb{Z}}_{0}^{\mathbb{Z}} }} F_{k}\Bigg(\sum \limits _{k \in {{\mathbb{Z}}_{0}^{\mathbb{Z}} }}y_{k}e^{i \langle k,\theta \rangle }\Bigg)e^{i \langle k,\theta \rangle }=\sum \limits _{k \in {{\mathbb{Z}}_{0}^{\mathbb{Z}} }} F_{k}(Y(\theta))e^{i \langle k,\theta \rangle}=\widetilde{F}(\theta,Y(\theta))
\end{align*}
is the shell function of $F(x,y(x))$.\qed

In the proof of the intersection property of the mapping $\mathfrak{M}$, we will use the generating function of the almost periodic monotonic twist map. Here, we describe its properties.

Consider an area-preserving almost periodic monotonic twist map defined in a horizontal strip $\mathbb{R}^1\times I\subseteq \mathbb{R}^2$ as follows :
\begin{equation}\label{g14}
\begin{array}{ll}
\phi : \mathbb{R}^1\times I\rightarrow \mathbb{R}^2 :\ \ \  \ \ \begin{array}{ll}
x_2 = x_1+\phi_1(x_1,y_1),\\[0.2cm]
y_2 =\phi_2(x_1,y_1),
\end{array}
\end{array}
\end{equation}
where $\phi_1(x_1,y_1),\phi_2(x_1,y_1)\in AP(\omega)$ as functions of $x_1$. Using the monotonic twist condition ${{\partial\phi_1}\over{\partial{y_1}}}(x_1,y_1)>0$, we invert the first equation of (\ref{g14}) and substitute it into the second equation of (\ref{g14}) to obtain
\begin{equation}\label{g15}
\begin{array}{ll}
\begin{array}{ll}
y_1 = \psi_1(x_2-x_1,x_1),\\[0.2cm]
y_2 =\psi_2(x_2-x_1,x_1),
\end{array}
\end{array}
\end{equation}
where $\psi_1$ is the inverse of $\phi_1$ and $\psi_2=\phi_2\big(x_1,\psi_1(x_2-x_1,x_1)\big)$.

By (\ref{g14}) and (\ref{g15}), we have
\begin{equation}\label{g16}
\begin{array}{ll}
\begin{array}{ll}
1={{\partial\phi_1}\over{\partial{y_1}}}{{\partial\psi_1}\over{\partial{x_2}}},\\[0.2cm]
1+{{\partial\phi_1}\over{\partial{x_1}}}+{{\partial\phi_1}\over{\partial{y_1}}}{{\partial\psi_1}\over{\partial{x_1}}}=0,\\[0.2cm]
{{\partial\psi_2}\over{\partial{x_1}}}={{\partial\phi_2}\over{\partial{x_1}}}+{{\partial\phi_2}\over{\partial{y_1}}}{{\partial\psi_1}\over{\partial{x_1}}}\\[0.2cm]
{{\partial\psi_2}\over{\partial{x_2}}}={{\partial\phi_2}\over{\partial{y_1}}}{{\partial\psi_1}\over{\partial{x_2}}}.
\end{array}
\end{array}
\end{equation}
Since ${{\partial\phi_1}\over{\partial{y_1}}}(x_1,y_1)>0$, then we have

\begin{eqnarray*}
{{\partial\psi_1}\over{\partial{x_2}}}+ {{\partial\psi_2}\over{\partial{x_1}}}&=& {{\partial\psi_1}\over{\partial{x_2}}}+{{\partial\phi_2}\over{\partial{x_1}}}+{{\partial\phi_2}\over{\partial{y_1}}}{{\partial\psi_1}\over{\partial{x_1}}}\\[0.2cm]
&=& {1\over {{{\partial\phi_1}\over{\partial{y_1}}}}}+{{\partial\phi_2}\over{\partial{x_1}}}+{{\partial\phi_2}\over{\partial{y_1}}}\Big(-1-{{\partial\phi_1}\over{\partial{x_1}}}\Big){1\over {{{\partial\phi_1}\over{\partial{y_1}}}}}\\[0.2cm]
&=& {1\over {{{\partial\phi_1}\over{\partial{y_1}}}}}\Bigg(1-\bigg[\Big(1+{{\partial\phi_1}\over{\partial{x_1}}}\Big){{\partial\phi_2}\over{\partial{y_1}}}-
{{\partial\phi_2}\over{\partial{x_1}}}{{{\partial\phi_1}\over{\partial{y_1}}}}\bigg]\Bigg).
\end{eqnarray*}
Using the area preservation property of $\phi$, we have
\begin{equation*}
\left|\begin{matrix}
1+{{\partial\phi_1}\over{\partial{x_1}}} & {{{\partial\phi_1}\over{\partial{y_1}}}}\\[0.2cm]
{{\partial\phi_2}\over{\partial{x_1}}} & {{{\partial\phi_2}\over{\partial{y_1}}}}
\end{matrix}\right|=\Big(1+{{\partial\phi_1}\over{\partial{x_1}}}\Big){{\partial\phi_2}\over{\partial{y_1}}}-
{{\partial\phi_2}\over{\partial{x_1}}}{{{\partial\phi_1}\over{\partial{y_1}}}}=1.
\end{equation*}
Hence,
\begin{equation}\label{g17}
{{\partial\psi_1}\over{\partial{x_2}}}+ {{\partial\psi_2}\over{\partial{x_1}}}=0.
\end{equation}
Therefore, there exists a generating function $h(x_1,x_2)$ defining the map implicitly
\begin{equation}\label{g18}
\begin{array}{ll}
\begin{array}{ll}
y_1 = -h_{x_1}(x_1,x_2),\\[0.2cm]
y_2 =h_{x_2}(x_1,x_2)
\end{array}
\end{array}
\end{equation}
with $h_{12}<0,$ defined in
$$\Omega=\Big\{(x_1,x_2)\in\mathbb{R}^2\ :\ a(x_1)<x_2<b(x_1)\Big\}$$
corresponding to the domain of definition of the map.

Integrating the vector field $(\psi_1,\psi_2)$ along a path from some point $(x_{01},x_{02})\in\Omega$ to $(x_{1},x_{2})\in\Omega$ and using the nonresonance condition (\ref{b6}), we obtain the generating function $h$ is of the form
\begin{equation}\label{g19}
h(x_1,x_2)=\beta x_1+H(x_2-x_1,x_1),
\end{equation}
where $H$ is real analytic almost periodic with the frequency  $\omega=(\cdots,\omega_\lambda,\cdots)$ in the second variable. When the map $\phi$ is exact, then $\beta=0$.

Now we are going to prove the following lemma.

\begin{lemma}\label{lem7.10}
If the mapping $\mathfrak{M}$ given by (\ref{M}) is an exact symplectic map, then it has intersection property.
\end{lemma}
\Proof Since the mapping $\mathfrak{M}$ given by (\ref{M}) is an exact symplectic map and it is also a monotonic twist map, from the above, there is a generating function $H$ such that the mapping $\mathfrak{M}$ can be rewritten by
\begin{equation}
y= -\frac{\partial }{\partial x} H(x_1-x, x),\quad y_1 = \frac{\partial}{\partial x_1}H(x_1-x,x),
\end{equation}
where $H$ is real analytic almost periodic with the frequency  $\omega=(\cdots,\omega_\lambda,\cdots)$ in the second variable.

Now we prove the intersection property of the mapping $\mathfrak{M}$, that is, given any continuous almost periodic curve $\Gamma : y=y(x)$, we need to prove that $\mathfrak{M}(\Gamma)\cap\Gamma\ne\emptyset.$  Define  two sets $\mathbb{B}$ and ${\mathbb{B}}_1$ : the set $\mathbb{B}$ is bounded by four curves $\big\{(x,y) : x=t\big\}$, $\big\{(x,y) : x=T\big\}$, $\big\{(x,y) : y=y_{*}\big\}$ and $\big\{(x,y) : y=y(x)\big\}$, the set ${\mathbb{B}}_1$ is bounded by four curves   $\big\{(x,y) : x=t\big\}$, $\big\{(x,y) : x=T\big\}$, $\big\{(x,y) : y=y_{*}\big\}$ and the image of $\Gamma$ under $\mathfrak{M}$. Here we choose $y_{*}<\min\limits _{t\leq x\leq T} y(x)$ and $y_{*}$ is smaller than the image of $\Gamma$ under $\mathfrak{M}$ as $t\leq x\leq T$. It is easy to show that the difference of the areas of ${\mathbb{B}}_1$ and $\mathbb{B}$ is
$$\Delta(t,T)=\int _{t}^{T}y_{1}d x_1-\int _{t}^{T}y d x=H(x_{1}(T)-T,T)-H(x_{1}(t)-t,t).$$
From the definition of $\mathfrak{M}$ and Lemma \ref{lem 7.1}, we know that $x_{1}(T)-T=y(T)+f(T,y(T))$ is almost periodic in $T$ and $x_{1}(t)-t=y(t)+f(t,y(t))$ is almost periodic in $t$. Hence using Lemma \ref{lem 7.1} again, it follows that  $\Delta(t,T)$ is almost periodic in $t$ and $T$.

Hence there are at least two pairs of $(t_1,T_1)$ and $(t_2,T_2)$ such that $\Delta(t_1,T_1)<0,\Delta(t_2,T_2)>0.$ The intersection property of $\mathfrak{M}$ follows from this fact, which proves the lemma. \qed

\begin{remark}
The intersection property is more flexible than exact symplecticity. For example, if a mapping $\mathcal{M}$ has the intersection property, then the transformed mapping of $\mathcal{M}$ under a diffeomorphism has also  the intersection property.
\end{remark}

Now we are in a position to  state our main result.
\begin{theorem}\label{thm2.11}
Suppose that the  almost periodic mapping $\mathfrak{M}$ given by (\ref{M}) has the intersection property, and for every $y$, $f(\cdot,y),g(\cdot,y)\in AP_r(\omega)$ with $\omega$ satisfying the nonresonance condition (\ref{b6}), and the corresponding shell functions $F(\theta,y),G(\theta,y)$ are real
analytic in the domain $D(r,s)=\big\{(\theta,y)\in \mathbb{C}^\mathbb{Z} \times \mathbb{C}\ :\ |\Im\,\theta|_{\infty}<r, |y-\alpha|<s\big\}$ with $\alpha$ satisfying (\ref{b7}). Then for each positive $\bar{\varepsilon}$, there is a positive $\varepsilon_0=\varepsilon_0(\bar{\varepsilon},r,s,m,\gamma,\Delta)$ such that if $f,g$ satisfy  the following smallness condition
$$\|f\|_{m,r,s}+\|g\|_{m,r,s}<{{\varepsilon}}_0,$$
then the almost periodic mapping $\mathfrak{M}$ has an invariant curve $\mathbf{\Gamma _{0}}$ with the form
$$
\begin{array}{ll}
\left\{\begin{array}{ll}
x=x'+\varphi(x'),\\[0.2cm]
y=\psi(x'),
 \end{array}\right.\ \ \
\end{array}
$$
where $\varphi, \psi$ are almost periodic  with the frequency $\omega=(\cdots,\omega_\lambda,\cdots)$, and the invariant curve $\mathbf{\Gamma _{0}}$ is of the form $y=\phi(x)$ with $\phi\in AP_{r'}(\omega)$ for some $r'<r$, and $\|\phi-\alpha\|_{m',r'}<\bar{\varepsilon},\ 0<m'<m$. Moreover,\ the  restriction of $\mathfrak{M}$ onto $\mathbf{\Gamma _{0}}$ is
 $$\mathfrak{M}|_{\mathbf{\Gamma _{0}}}: x_{1}^\prime=x^\prime +\alpha.$$
\end{theorem}

\begin{remark}
If all the conditions of Theorem \ref{thm2.11} hold,\ then the mapping $\mathfrak{M}$ has many invariant curves $\mathbf{\Gamma _{0}}$ which can be labeled by the form $$\mathfrak{M}|_{\mathbf{\Gamma _{0}}}: x_{1}^\prime = x' +\alpha$$
 of the restriction of $\mathfrak{M}$ onto $\mathbf{\Gamma _{0}}.$\ In fact,\ given any $\alpha$ satisfying the inequalities (\ref{b7}), there exists an invariant curve $\mathbf{\Gamma _{0}}$ of $\mathfrak{M}$ which is almost periodic  with the frequency  $\omega=(\cdots,\omega_\lambda,\cdots)$,\ and the restriction of $\mathfrak{M}$ onto $\mathbf{\Gamma _{0}}$ has the form
$$\mathfrak{M}|_{\mathbf{\Gamma _{0}}}: x_{1}^\prime=x' +\alpha.$$
\end{remark}

\section{The measure estimate}
In this section the measure estimate of the rotation numbers $\alpha$ satisfying inequalities  (\ref{b7}) will be given. Firstly, we give some useful  lemmas.
\begin{lemma}[Lemma\ 2\ in {\cite{Poschel90}}]\label{lem4.1}
For every given approximation function $\Theta$, there exists an approximation
function $\Delta$ such that
$$\sum \limits_{A\in\mathcal{S},|A|=n} {1\over {\Delta([A])}}\leq {{2N_0}\over {\Theta(t_n)}},\ \ \ n\geq 1,$$
where $N_0, t_n$ are given in Lemma \ref{lem2.11}.
\end{lemma}

\begin{lemma}[Lemma\ 4\ in {\cite{Poschel90}}]\label{lem4.3}
There is an approximation function $\Delta$ such that
$$\sum \limits_{l\in\mathbb{Z}^n} {1\over {\Delta(|l|)}}\leq {K}^{n\log\log n}$$
for all sufficiently large $n$ with some constant ${K}$, where $l=(l_1,l_2,\cdots,l_n)$ and $|l|=|l_1|+|l_2|+\cdots+|l_n|$.
\end{lemma}

\begin{remark}
From the proofs of Lemma \ref{lem4.1} and  Lemma \ref{lem4.3}, we know that there is the same approximation function $\Delta$ such that Lemma \ref{lem4.1} and  Lemma \ref{lem4.3} hold simultaneously. The detail proofs of Lemma \ref{lem4.1} and  Lemma \ref{lem4.3} can be found in {\cite{Poschel90}}.
\end{remark}

\begin{theorem}\label{thm2.12}
There is an approximation function $\Delta$ such that for suitable $\gamma$,\ the set of $\alpha$ satisfying (\ref{b7}) has positive measure.
\end{theorem}

\Proof Choose the frequency  $\omega=(\cdots,\omega_\lambda,\cdots)$  satisfying the nonresonance condition (\ref{b6}) and let $\mathcal{D}_{\gamma}^{\omega}$ denote the set of all $\alpha\in\mathbb{R}$ satisfying (\ref{b7}) with the fixed $\gamma$. Then $\mathcal{D}_{\gamma}^{\omega}$ is the complement of the open dense set $\mathcal{R}_{\gamma}^{\omega}$, where
\begin{eqnarray*}
\mathcal{R}_{\gamma}^{\omega}&=&\bigcup \limits _{\substack{  k \in \mathbb{Z}_{0}^{\mathbb{Z}}\backslash\{0\}\\  j \in \mathbb{Z}}}\mathcal{R}_{\omega,\gamma}^{k,j}\\&=&\bigcup \limits _{\substack{  k \in \mathbb{Z}_{0}^{\mathbb{Z}}\backslash\{0\}\\  j \in \mathbb{Z}}}\Big\{\alpha\in [a+\gamma,b-\gamma]:\big|{\langle k,\omega \rangle {\alpha \over {2\pi}}-j}\big|<{\gamma \over {\Delta([[k]])\Delta(|k|)}}\Big\}.
\end{eqnarray*}

Now we estimate the measure of the set $\mathcal{R}_{\omega,\gamma}^{k,j}$. Since $k \in \mathbb{Z}_{0}^{{\mathbb{Z}}}\backslash\{0\}$, set $|k_{\max}|=\max\limits_{\lambda\in \mbox{supp}\, k}|k_\lambda|$, then there exists some  $m\in \mbox{supp}\,k$ such that $|k_m|=|k_{\max}|$, and
$1\leq{{|k|}\over {|k_{\max}|}}\leq |\mbox{supp}\,k|.$  Therefore, we have
\begin{eqnarray*}
\mathcal{R}_{\omega,\gamma}^{k,j}&=&\Big\{\alpha\in [a+\gamma,b-\gamma]:\big|{\langle k,\omega \rangle {\alpha \over {2\pi}}-j}\big|<{\gamma \over {\Delta([[k]])\Delta(|k|)}}\Big\}\\
&=& \Big\{\alpha\in [a+\gamma,b-\gamma]:\big|k_{\max}\omega_m{\alpha \over {2\pi}}+\sum_{\lambda\not =m} k_\lambda\omega_\lambda{\alpha \over {2\pi}}-j\big|<{\gamma \over {\Delta([[k]])\Delta(|k|)}}\Big\}\\
&=& \Big\{\alpha\in [a+\gamma,b-\gamma]:|k_{\max}||\omega_m||\alpha+b_{j}|<{2\pi\gamma \over {\Delta([[k]])\Delta(|k|)}}\Big\}\\
&=& \Big\{\alpha\in [a+\gamma,b-\gamma]:-b_{j}-\delta_k<\alpha<-b_{j}+\delta_k\Big\},
\end{eqnarray*}
where $b_{j}={1\over {k_{\max}\omega_m}}\Big\{\sum\limits_{\lambda\not =m} k_\lambda\omega_\lambda\alpha-2\pi j\Big\}$ and
$
\delta_k = {2\pi\gamma \over {\Delta([[k]])\Delta(|k|)}}\,{1\over {|k_{\max}||\omega_m|}}.
$
Hence,
$$\mbox{meas}\big(\mathcal{R}_{\omega,\gamma}^{k,j}\big)\le 2\delta_k = {4\pi\gamma \over {\Delta([[k]])\Delta(|k|)}}\,{1\over {|k_{\max}||\omega_m|}}={4\pi\gamma \over {|k|\Delta([[k]])\Delta(|k|)}}\,{|k|\over {|k_{\max}|}}\,{1\over{|\omega_m|}}.
$$
Since $1\leq{{|k|}\over {k_{\max}}}\leq |\mbox{supp}\,k|,$ then we have the following measure estimate
$$\gamma^{-1}\mbox{meas}\big(\mathcal{R}_{\omega,\gamma}^{k,j}\big)\leq {C_0 \over {|k|\Delta([[k]])\Delta(|k|)}}$$
with some positive constant $C_0$.

Next we estimate the measure of the set $\mathcal{R}_{\gamma}^{\omega}$. Since for $\alpha\in\mathcal{R}_{\omega,\gamma}^{k,j}$,
$$\big|{\langle k,\omega \rangle {\alpha \over {2\pi}}-j}\big|<{\gamma \over {\Delta([[k]])\Delta(|k|)}},$$
then we have
$$|j|\leq \big|\langle k,\omega \rangle \big|{{\alpha}\over {2\pi}}+{\gamma \over {\Delta([[k]])\Delta(|k|)}}\leq c_0|k|,$$
where $c_0$ is a constant independent of $k$. Thus
\begin{eqnarray*}
\gamma^{-1}\mbox{meas}(\mathcal{R}_{\gamma}^{\omega}) &\leq & \sum\limits_{k \in \mathbb{Z}_{0}^{\mathbb{Z}}\backslash\{0\}}\sum\limits_{\substack{j\in \mathbb{Z}\\ |j|\leq c_0|k|}}\gamma^{-1}\mbox{meas}\big(\mathcal{R}_{\omega,\gamma}^{k,j}\big)\\[0.2cm]
&\leq& \sum\limits_{k \in \mathbb{Z}_{0}^{\mathbb{Z}}\backslash\{0\}}\sum\limits_{\substack{j\in \mathbb{Z}\\ |j|\leq c_0|k|}}{C_0 \over {|k|\Delta([[k]])\Delta(|k|)}}\leq C_1\sum\limits_{k \in  \mathbb{Z}_{0}^{\mathbb{Z}}\backslash\{0\}}{1 \over {\Delta([[k]])\Delta(|k|)}}\\[0.0001cm]
&\leq& C_1\sum\limits_{A\in\mathcal{S}}\Bigg({1\over \Delta([A])}\sum\limits_{\substack{\mbox{supp}\,k\subseteq A\\ k\neq 0}}{1\over {\Delta(|k|)}}\Bigg)\\[0.0001cm]
&\leq& C_1\sum\limits_{n=1}^{+\infty}\Bigg\{\Bigg(\sum\limits_{A\in\mathcal{S},|A|=n}{1\over \Delta([A])}\Bigg)\sum\limits_{l\in \mathbb{Z}^n\backslash\{0\}}{1\over {\Delta(|l|)}}\Bigg\}
\end{eqnarray*}
with some positive constant $C_1$. Thus the sum is broken up with respect to the cardinality and the weight of the spatial
components of $\mathcal{S}$. Each of these factors is now studied separately.
By applying Lemma \ref{lem4.1}, we have
$$\sum \limits_{A\in\mathcal{S},|A|=n} {1\over {\Delta([A])}}\leq {{2N_0}\over {\Theta(t_n)}},\ \ \ n\geq 1.$$
From Lemma \ref{lem4.3}, we get
$$\sum \limits_{l\in\mathbb{Z}^n\backslash\{0\}} {1\over {\Delta(|l|)}}\leq {K}^{n\log\log n}.$$
Summarizing all our estimates so far we arrive at
\begin{eqnarray*}
\gamma^{-1}\mbox{meas}(\mathcal{R}_{\gamma}^{\omega}) &\leq& C_1\sum\limits_{n=1}^{+\infty}\Bigg\{\Bigg(\sum\limits_{A\in\mathcal{S},|A|=n}{1\over \Delta([A])}\Bigg)\sum\limits_{l\in \mathbb{Z}^n\backslash\{0\}}{1\over {\Delta(|l|)}}\Bigg\}\\[0.2cm]
&\leq& C+C\sum\limits_{n=n_0}^{+\infty}{{{K}^{n\log\log n}}\over{\Theta(t_n)}}
\end{eqnarray*}
with some constant $C$ and $n_0$ so large that $t_n\geq n\log^{\varrho-1} n$ for $n\geq n_0, \varrho>2$ by hypotheses.
Here we are still free to choose a suitable approximation function $\Theta$, and choose
$$\Theta(t)=\exp\Bigg({t\over{\log\,t\log^{\varrho-1}\log\,t}}\Bigg),\ \ \ t>e, \varrho>2,$$
the infinite sum does converge. Thus there is an approximation function $\Delta$ such that
$$\gamma^{-1}\mbox{meas}(\mathcal{R}_{\gamma}^{\omega})<+\infty.$$
Hence,
$$\mbox{meas}(\mathcal{R}_{\gamma,\tau}^{\omega})\leq O(\gamma)$$
and
$$\mbox{meas}(\mathcal{D}_{\gamma,\tau}^{\omega})\to b-a \ \ \ \ \ \mbox{as} \ \ \ \ \ \gamma\to 0.$$
This completes the proof of Theorem \ref{thm2.12}.\qed

\section{The KAM step}
The proof of Theorem \ref{thm2.11} is based on the KAM approach, is to find a sequence of changes
of variables such that the transformed mapping of $\mathfrak{M}$ will be closer to
$$\begin{array}{ll}
\left\{\begin{array}{ll}
x_1=x+y,\\[0.2cm]
y_1=y
 \end{array}\right.
\end{array}$$
than the previous one in the narrower domain. In the following, we will give a construction
of such transformation.

\subsection{Constructions of $u,v$}

In this subsection, we shall construct a change of variables $\mathfrak{U}$
\begin{equation}\label{c1}
\begin{array}{ll}
\left\{\begin{array}{ll}
x = \xi+u(\xi,\eta),\\[0.2cm]
y = \eta+v(\xi,\eta),
\end{array}\right.
\end{array}
\end{equation}
where $u$ and $v$ are real analytic  and  almost periodic in $\xi$. Under this transformation, the original mapping $\mathfrak{M}$ is changed into the form
\begin{equation}\label{c2}
\begin{array}{ll}
\mathfrak{{U}}^{-1} \mathfrak{M}\, \mathfrak{U}\ :\ \left\{\begin{array}{ll}
\xi_1=\xi+\eta+f_{+}(\xi,\eta),\\[0.2cm]
\eta_1=\eta+g_{+}(\xi,\eta),
 \end{array}\right.
\end{array}
\end{equation}
where the functions $f_{+}$ and $g_{+}$ are real analytic  almost periodic functions in $\xi$ defined in a smaller domain $D(r_{+},s_{+})$ and $\|f_{+}\|_{m_+,r_+,s_+}+\|g_{+}\|_{m_+,r_+,s_+}$ is smaller than $\|f\|_{m,r,s}+\|g\|_{m,r,s}.$

In the following, $c_1,c_2,\cdots$ are positive constants depending on $\gamma,\Delta$ only. We
also assume that for each fixed $y$, $f(\cdot,y),g(\cdot,y)\in AP_r(\omega)$ and $f,g$ are real analytic in the domain $D(r,s)$ with $0<r<1,0<s<{1\over 2}$. Moreover, we assume that  $\omega=(\cdots,\omega_\lambda,\cdots),\alpha$ satisfy the nonresonance conditions (\ref{b6}) and (\ref{b7}). Let
$$\varepsilon=\|f\|_{m,r,s}+\|g\|_{m,r,s}.$$

We try to motivate the following constructions for $u,v$ first. From (\ref{c1}) and (\ref{c2}), it follows that
\begin{equation}\label{c4}
\left\{\begin{array}{ll}
f_{+}(\xi, \eta) = f(\xi + u, \eta + v) + v(\xi, \eta) + u(\xi, \eta) - u(\xi + \eta + f_{+},\eta+ g_{+}),\\[0.2cm]
g_{+}(\xi, \eta) = g(\xi + u, \eta + v) + v(\xi, \eta) - v(\xi + \eta + f_{+},\eta+ g_{+}),
\end{array}\right.
\end{equation}
which serve to define $f_{+},g_{+}$ implicity in $D(r_{+},s_{+})$.

In the following, we will determine the unknown functions $u$ and $v$.\ As one did in the periodic case, we may solve $u$ and $v$ from the following
equations
\begin{eqnarray*}
&&u(\xi + \alpha, \eta)-u(\xi, \eta) = v(\xi, \eta) + f(\xi,\eta), \\[0.2cm]
&&v(\xi + \alpha, \eta)-v(\xi, \eta) = g(\xi, \eta).
\end{eqnarray*}
These difference equation will introduce the small divisors.\ We can solve the second equation only if the mean value of $g$ over  the first variable vanishes.

For this reason we define $u,v$ as the solutions of the following modified homological equations
\begin{equation}\label{c5}
\left\{\begin{array}{ll}
u(\xi + \alpha, \eta)-u(\xi, \eta) = v(\xi, \eta) + f(\xi,\eta),\\[0.2cm]
v(\xi + \alpha, \eta)-v(\xi, \eta) = g\xi, \eta)-[g](\eta).
\end{array}\right.
\end{equation}
Here [\ ] denotes the mean value of a function over the first variable.\

\subsection{Estimates of $u,v$}

In order to solve $u,v$ from the equations (\ref{c5}), we will meet the following difference equation (the so-called
homological equation) :
\begin{equation}\label{b8}
l(x+\alpha)-l(x)=h(x),
\end{equation}
where $h\in AP_r(\omega)$. Let us first study this equation.

\begin{lemma}\label{lem2.13}
Suppose that $h\in AP_r(\omega)$ and $\omega=(\cdots,\omega_\lambda,\cdots),\alpha$ satisfy the nonresonance condition
\begin{equation}\label{b9}
\Big|{\langle k,\omega \rangle {\alpha \over {2\pi}}-j}\Big|\geq {\gamma \over {\Delta([[k]])\Delta(|k|)}},\ \ \ \  \mbox{for all}\ \  k \in \mathbb{Z}_{0}^{{\mathbb{Z}}}\backslash\{0\},\ \  j \in \mathbb{Z},\ \ \gamma>0.
\end{equation}
Then for any $0<r'<r$, the difference equation (\ref{b8}) has the unique
solution $l\in AP_{r'}(\omega)$ with $\lim \limits_{T\rightarrow\infty}{1\over T}\int_{0}^{T} l(x)dx=0$ if and only if
\begin{equation}\label{b10}
\lim \limits_{T\rightarrow\infty}{1\over T}\int_{0}^{T} h(x)dx=0.
\end{equation}
In this case, we have the following estimate
\begin{equation}\label{b11}
\|l\|_{m',r'}\leq {\gamma^{-1}}\,\Gamma_{0}(r-r')\,\Gamma_{0}(m-m')\|h\|_{m,r}
\end{equation}
for $0<m'<m,\ \Gamma_{0}(\rho)=\sup \limits_{t\geq 0}\Delta(t)e^{-\rho t}.$
\end{lemma}
\noindent\textbf{Proof}: From $h\in AP_{r}(\omega)$ and (\ref{b10}), we know that $h$ can be represented by
\begin{equation*}
\begin{array}{ll}
h(x)=\sum \limits _{A\in\mathcal{S}} h_{A}(x),
\end{array}
\end{equation*}
where
$$h_{A}(x)={\sum \limits_{\substack{k\neq 0\\ \mbox{supp}\, k \subseteq A}}} h_{k}\,e^{i \langle k,\omega \rangle x}.$$

Let
$$l(x)={\sum \limits_{A\in \mathcal{S}}}\ l_{A}(x),$$
where
$${ l_{A}(x)=\sum \limits_{ \mbox{supp}\, k \subseteq A}} l_{k}e^{i \langle k,\omega \rangle x}.$$
After straightforward calculations we obtain the relation between Fourier coefficients $h_{k}$ and $l_{k}$ as follows
$$
l_{k}={{h_{k}}\over {e^{i\langle k, \omega\rangle \alpha}-1}},\ \ \ k\neq 0,
$$
then $l$ is of the form
$$l(x)={\sum \limits_{A\in \mathcal{S}}}l_{A}(x)={\sum \limits_{A\in \mathcal{S}}}\ {\sum \limits_{\substack{k\neq 0\\ \mbox{supp}\, k \subseteq A}}} {{h_{k}}\over {e^{i\langle k, \omega\rangle \alpha}-1}}e^{i \langle k,\omega \rangle x},$$
which is the  uniquely determined Fourier expansion of the  wanted solution $l$  satisfying $l \in AP(\omega)$ with
$\lim \limits_{T\rightarrow\infty}{1\over T}\int_{0}^{T} l(x)dx=0.$

From (\ref{b9}), it follows that
$$\Big|e^{i\langle k,\omega\rangle \alpha}-1\Big|\geq \Big|{\langle k,\omega \rangle {\alpha \over {2\pi}}-j}\Big|\geq {\gamma \over {\Delta([[k]])\Delta(|k|)}}\ \ \ \  \mbox{for all}\ \  k \in \mathbb{Z}_{0}^{{\mathbb{Z}}}\backslash\{0\}.$$
Hence,
\begin{eqnarray*}
|l_A|_{r'} &\leq& {\sum \limits_{\substack{k\neq 0\\ \mbox{supp}\, k \subseteq A}}}{\gamma^{-1}}\Delta([[k]])\,\Delta(|k|)\,|h_{k}|\,e^{|k|r'}\\[0.001cm]
&\leq& {\sum \limits_{\substack{k\neq 0\\ \mbox{supp}\, k \subseteq A}}}{\gamma^{-1}}\,\Delta([[k]])\,\Delta(|k|)\,e^{-(r-r')|k|}\,|h_{k}|\,e^{|k|r}\\[0.001cm]
&\leq& {\gamma^{-1}}\Delta([A])\,\Gamma_{0}(r-r')\,|h_{A}|_{r},
\end{eqnarray*}
where $\Gamma_{0}(r-r')=\sup \limits_{t\geq 0}\Delta(t)e^{-(r-r')t}.$ Putting the spatial components together yields that
\begin{eqnarray*}
\|l\|_{m',r'} &\leq& {\gamma^{-1}} {\sum \limits_{A\in \mathcal{S}}}\Delta([A])\,\Gamma_{0}(r-r')\,|h_{A}|_{r}\,e^{m'[A]}\\[0.001cm]
&\leq& {\gamma^{-1}}\,\Gamma_{0}(r-r')\ {\sum \limits_{A\in \mathcal{S}}}\Delta([A])\,e^{-(m-m')[A]}\,|h_{A}|_{r}\,e^{m [A]}\\[0.001cm]
&\leq& {\gamma^{-1}}\,\Gamma_{0}(r-r')\,\Gamma_{0}(m-m')\|h\|_{m,r}
\end{eqnarray*}
for $0<m'<m,\ 0<r'<r,$
which completes the proof of the lemma.\qed

Now we can solve the functions $u$ and $v$ from (\ref{c5}) and give the estimates of them. In the first equation of (\ref{c5}) the  mean value over  the first variable must vanish on both sides. Hence we get the condition
\begin{equation}\label{c6}
[v](\eta)= -[f](\eta)
\end{equation}
for  the  mean value of $v$ over  the first variable. As a consequence,  we have
\begin{equation*}
\|[v]\|_{m,r,s}\leq \|f\|_{m,r,s}.
\end{equation*}
Lemma \ref{lem2.13} gives a unique solution $v(\xi,\eta)=\tilde{v}(\xi,\eta)$ of the second equation of (\ref{c5}) with $[\tilde{v}](\eta)=0$. This solution has the estimate
$$\|\tilde{v}\|_{m-\nu,r-\delta,s}\leq {\gamma^{-1}}\,\Gamma_{0}(\delta)\,\Gamma_{0}(\nu)\|g\|_{m,r,s}$$
for $0<\nu<m,\ 0<\delta<r.$ Define $v(\xi,\eta)=\tilde{v}(\xi,\eta)+[v](\eta)$, we obtain the  uniquely determined solution $v(\xi,\eta)$ of the second equation of (\ref{c5}).

Define $p(\xi,\eta)= \tilde{v}(\xi,\eta)+f(\xi,\eta)$, note that  $\tilde{v}(\xi,\eta)$ is defined in $D(r-\delta,s)$,\ then  $p(\xi,\eta)$ is well defined in $D(r-\delta,s)$. As a consequence we have
\begin{eqnarray*}
\|p\|_{m-\nu,r-\delta,s} &=&\|\tilde{v}+f\|_{m-\nu,r-\delta,s}\\[0.2cm]
&\leq & c_{1}\Gamma_{0}(\delta)\,\Gamma_{0}(\nu)(\|f\|_{m,r,s}+\|g\|_{m,r,s})
\end{eqnarray*}
and
\begin{eqnarray*}
p(\xi,\eta)-[p](\eta)
&=&\tilde{v}(\xi,\eta)+f(\xi,\eta)-[ \tilde{v}+f](\eta) \\[0.2cm]
&=&  \tilde{v}(\xi,\eta)+f(\xi,\eta)-[\tilde{v}](\eta)-[f](\eta) \\[0.2cm]
&=&  \tilde{v}(\xi,\eta)+f(\xi,\eta)+[v](\eta)\\[0.2cm]
&= &  v(\xi,\eta)+f(\xi,\eta).
\end{eqnarray*}
Hence, the first equation of (\ref{c5}) can  be written in  the  form
\begin{equation}\label{c7}
u(\xi + \alpha, \eta)-u(\xi, \eta) = p(\xi, \eta) - [p](\eta).
\end{equation}
Thus Lemma  \ref{lem2.13} gives a uniquely determined solution  $u$ of (\ref{c7}) with $[u]=0.$ For an estimate of $u$ we apply Lemma  \ref{lem2.13} to  (\ref{c7}) restricted  to $D(r-\delta,s)$ such that  in  (\ref{b11}) we have to  replace $h$ by $p$ and $r$ by $r-\delta.$

From the above discussions, we have
\begin{eqnarray}\label{c8}
\|u\|_{m-2\nu,r-2\delta,s}+\|v\|_{m-2\nu,r-2\delta,s} &\leq& c_{2}\Gamma_{0}^2(\delta)\,\Gamma_{0}^2(\nu)(\|f\|_{m,r,s}+\|g\|_{m,r,s})\nonumber\\[0.2cm]
&=& c_{2}\Gamma_{0}^2(\delta)\,\Gamma_{0}^2(\nu) \varepsilon
\end{eqnarray}
and by Cauchy's estimate
\begin{equation}\label{c9}
\left\|{{\partial u}\over {\partial \xi}}\right\|_{m-2\nu,r-3\delta,s}+\left\|{{\partial v}\over {\partial \xi}}\right\|_{m-2\nu,r-3\delta,s}\leq c_{2}\Gamma_{0}^2(\delta)\,\Gamma_{0}^2(\nu)\,{\varepsilon\over \delta},
\end{equation}
\begin{equation}\label{c28}
\left\|{{\partial u}\over {\partial \eta}}\right\|_{m-2\nu,r-2\delta,s-\rho}+\left\|{{\partial v}\over {\partial \eta}}\right\|_{m-2\nu,r-2\delta,s-\rho}
 \leq c_{2}\Gamma_{0}^2(\delta)\,\Gamma_{0}^2(\nu)\,{\varepsilon\over \rho}
\end{equation}
for $0<2\nu<m,0<3\delta<r$ and $0<\rho<s$.

For $0<m_+<m,0<r_+<r<1$ and $0<s_+<s<{1\over 2}$, let
$$\nu={1\over 10}(m-m_+),\ \ \ \ \delta={1\over 10}(r-r_+),\ \ \ \ \rho={1\over 10}(s-s_+).$$
Introduce some domains $D_j$ between $D(r_+,s_+)$ and $D(r,s)$ by
$$D_j=D(r-j\delta,s-j\rho)\ \ \ \ \mbox{for}\ \ 0\leq j\leq 10.$$
From (\ref{c8}) and (\ref{c9}), it follows that
\begin{eqnarray}\label{c10}
\mathfrak{U}^{-1}(D_{j+1})\subset D_j,\ \ \ \ \mathfrak{M}(D_{j+1})\subset D_j,\ \ \ \ \mathfrak{U}(D_{j+1})\subset D_j,
\end{eqnarray}
if
\begin{eqnarray}\label{c11}
\varepsilon^{-1}>c_{2}\Gamma_{0}^2(\delta)\,\Gamma_{0}^2(\nu)\,\max\Big\{{1\over \delta},{1\over \rho}\Big\}.
\end{eqnarray}
Hence the mapping $\mathfrak{U}^{-1}\mathfrak{M}\,\mathfrak{U}$ is well defined in the domain $D(r_+,s_+)$ and maps this domain into $D_6$.

\subsection{The estimates on the new perturbations}

In this subsection the estimates on the new perturbations $f_{+},g_{+}$ are given.

Firstly, if $u(\cdot,\eta), v(\cdot,\eta) \in AP(\omega)$, then one can prove that $f_{+}(\cdot,\eta), g_{+}(\cdot,\eta)\in AP(\omega)$ are well defined by (\ref{c4}). Indeed, since $f(\cdot,y),g(\cdot,y),u(\cdot,\eta), v(\cdot,\eta)\in AP(\omega)$, then we have
$$f(x,y)=\sum \limits _{k \in {{\mathbb{Z}}_{0}^{\mathbb{Z}} }}f_{k}(y)e^{i \langle k,\omega \rangle x},\  F(\theta,y)=\sum \limits _{k \in {{\mathbb{Z}}_{0}^{\mathbb{Z}} }} f_{k}(y)e^{i \langle k,\theta \rangle },$$
$$g(x,y)=\sum \limits _{k \in {{\mathbb{Z}}_{0}^{\mathbb{Z}} }} g_{k}(y)e^{i \langle k,\omega \rangle t},\ G(\theta,y)=\sum \limits _{k \in {{\mathbb{Z}}_{0}^{\mathbb{Z}} }} g_{k}(y)e^{i \langle k,\theta \rangle },$$
and
$$u(\xi,\eta)=\sum \limits _{k \in {{\mathbb{Z}}_{0}^{\mathbb{Z}} }}u_{k}(\eta)e^{i \langle k,\omega \rangle \xi},\  U(\theta,\eta)=\sum \limits _{k \in {{\mathbb{Z}}_{0}^{\mathbb{Z}} }} u_{k}(\eta)e^{i \langle k,\theta \rangle },$$
$$v(\xi,\eta)=\sum \limits _{k \in {{\mathbb{Z}}_{0}^{\mathbb{Z}} }} v_{k}(\eta)e^{i \langle k,\omega \rangle \xi},\ V(\theta,\eta)=\sum \limits _{k \in {{\mathbb{Z}}_{0}^{\mathbb{Z}} }} v_{k}(\eta)e^{i \langle k,\theta \rangle }.$$
Then
$$f(\xi + u, \eta + v)=\sum \limits _{k \in {{\mathbb{Z}}_{0}^{\mathbb{Z}} }}f_{k}\Bigg(\eta+\sum \limits _{k \in {{\mathbb{Z}}_{0}^{\mathbb{Z}} }} v_{k}(\eta)e^{i \langle k,\omega \rangle \xi}\Bigg)e^{i \langle k,\omega \rangle \Big(\xi+\sum \limits _{k \in {{\mathbb{Z}}_{0}^{\mathbb{Z}} }}u_{k}(\eta)e^{i \langle k,\omega \rangle \xi}\Big)},$$
$$g(\xi + u, \eta + v)=\sum \limits _{k \in {{\mathbb{Z}}_{0}^{\mathbb{Z}} }}g_{k}\Bigg(\eta+\sum \limits _{k \in {{\mathbb{Z}}_{0}^{\mathbb{Z}} }} v_{k}(\eta)e^{i \langle k,\omega \rangle \xi}\Bigg)e^{i \langle k,\omega \rangle \Big(\xi+\sum \limits _{k \in {{\mathbb{Z}}_{0}^{\mathbb{Z}} }}u_{k}(\eta)e^{i \langle k,\omega \rangle \xi}\Big)}.$$
Hence,
\begin{align*}
&\sum \limits _{k \in {{\mathbb{Z}}_{0}^{\mathbb{Z}} }}f_{k}\Bigg(\eta+\sum \limits _{k \in {{\mathbb{Z}}_{0}^{\mathbb{Z}} }} v_{k}(\eta)e^{i \langle k,\theta \rangle}\Bigg)e^{i \Big\langle k,\theta+\omega\sum \limits _{k \in {{\mathbb{Z}}_{0}^{\mathbb{Z}} }}u_{k}(\eta)e^{i \langle k,\theta \rangle }\Big\rangle }\\[0.2cm]
&=\sum \limits _{k \in {{\mathbb{Z}}_{0}^{\mathbb{Z}} }}f_{k}\Big(\eta+V(\theta,\eta)\Big)e^{i \langle k,\theta+\omega U(\theta,\eta) \rangle } =F(\theta+\omega U(\theta,\eta), \eta+V(\theta,\eta)),\\[0.2cm]
&\sum \limits _{k \in {{\mathbb{Z}}_{0}^{\mathbb{Z}} }}g_{k}\Bigg(\eta+\sum \limits _{k \in {{\mathbb{Z}}_{0}^{\mathbb{Z}} }} v_{k}(\eta)e^{i \langle k,\theta \rangle}\Bigg)e^{i \Big\langle k,\theta+\omega\sum \limits _{k \in {{\mathbb{Z}}_{0}^{\mathbb{Z}} }}u_{k}(\eta)e^{i \langle k,\theta \rangle }\Big\rangle }\\[0.2cm]
&=\sum \limits _{k \in {{\mathbb{Z}}_{0}^{\mathbb{Z}} }}g_{k}\Big(\eta+V(\theta,\eta)\Big)e^{i \langle k,\theta+\omega U(\theta,\eta)\rangle }=G(\theta+\omega U(\theta,\eta), \eta+V(\theta,\eta))
\end{align*}
are the shell functions of $f(\xi + u, \eta + v),g(\xi + u, \eta + v)$, respectively. That is to say, $f(\xi + u, \eta + v),g(\xi + u, \eta + v)$ are real analytic  almost periodic functions in $\xi$. From Lemma \ref{lem2.4}, we know $f(\xi + u, \eta + v) + v(\xi, \eta) + u(\xi, \eta), g(\xi + u, \eta + v) + v(\xi, \eta)$ are also real analytic  almost periodic functions in $\xi$. Denote $\tilde{\phi}(\xi,\eta)=f(\xi + u, \eta + v) + v(\xi, \eta) + u(\xi, \eta),\tilde{\psi}(\xi,\eta)=g(\xi + u, \eta + v) + v(\xi, \eta)$, and the corresponding shell functions are $\tilde{\Phi}(\theta,\eta),\tilde{\Psi}(\theta,\eta)$ which are  period $2\pi$ in each of the variables $\theta_{\lambda}\,( \lambda\in \mathbb{Z} )$. If the unknown functions $f_+(\xi,\eta), g_+(\xi,\eta)$ are represented by $F_+(\theta,\eta), G_+(\theta,\eta)$, where $f_+(\xi,\eta)=F_+(\omega\xi,\eta),g_+(\xi,\eta)=G_+(\omega\xi,\eta)$, the conditions for $F_+,G_+$ become
\begin{equation}\label{c30}
\left\{\begin{array}{ll}
F_+(\theta,\eta) = \tilde{\Phi}(\theta,\eta) - U(\theta + \omega\eta + \omega F_+(\theta,\eta),\eta+ G_+(\theta,\eta)),\\[0.2cm]
G_+(\theta,\eta)) = \tilde{\Psi}(\theta,\eta) - V(\theta + \omega \eta + \omega F_+(\theta,\eta),\eta+G_+(\theta,\eta)).
\end{array}\right.
\end{equation}
By (\ref{c9}), (\ref{c28}), (\ref{c11}) and  the implicit function theorem, we know $F_+(\theta,\eta),\\ G_+(\theta,\eta)$ are well defined by (\ref{c30}), and have period $2\pi$ in each of the variables $\theta_{\lambda}\,( \lambda\in \mathbb{Z} )$.

Hence, from the constructions of $u,v$ and the relationships between $f_{+},g_{+}$ and $F_+,G_+$, we know that the new perturbations $f_{+}(\cdot,\eta), g_{+}(\cdot,\eta)\in AP(\omega)$. Similar to the proof in \cite[chapter 3]{Siegel97}, one may prove that $f_+$ and $g_+$  are real analytic in $D(r_+,s_+)$. Moreover, we have
$$\|f_+\|_{m_+,r_+,s_+}<4\delta,\ \ \ \ \|g_+\|_{m_+,r_+,s_+}<4\rho.$$

In the following, we will prove that $\|f_+\|_{m_+,r_+,s_+}+\|g_+\|_{m_+,r_+,s_+}$ is much smaller than $\|f\|_{m,r,s}+\|g\|_{m,r,s}$.

From (\ref{c4}) and (\ref{c5}), we have
\begin{equation}\label{c22}
\left\{\begin{array}{ll}
f_{+}(\xi, \eta) = u(\xi+\alpha, \eta) - u(\xi_1,\eta_1)+f(\xi + u, \eta + v) -f(\xi,\eta) ,\\[0.3cm]
g_{+}(\xi, \eta) = v(\xi+\alpha, \eta)- v(\xi_1,\eta_1)+g(\xi + u, \eta + v)-g(\xi,\eta)+[g](\eta).
\end{array}\right.
\end{equation}

We obtain from (\ref{c8}), (\ref{c9}), (\ref{c28}), the estimates
\begin{eqnarray*}
\|u\|_{m_+,r_+,s_+}+\|v\|_{m_+,r_+,s_+} < c_{2}\Gamma_{0}^2\Big({{r-r_+}\over 10}\Big)\,\Gamma_{0}^2\Big({{m-m_+}\over 10}\Big) \varepsilon,
\end{eqnarray*}
\begin{equation*}
\left\|{{\partial u}\over {\partial \xi}}\right\|_{m_+,r_+,s_+}+\left\|{{\partial v}\over {\partial \xi}}\right\|_{m_+,r_+,s_+}< c_{3}\Gamma_{0}^2\Big({{r-r_+}\over 10}\Big)\,\Gamma_{0}^2\Big({{m-m_+}\over 10}\Big)\,{\varepsilon\over {r-r_+}},
\end{equation*}
\begin{align*}
\left\|{{\partial u}\over {\partial \eta}}\right\|_{m_+,r_+,s_+}+\left\|{{\partial v}\over {\partial \eta}}\right\|_{m_+,r_+,s_+}
& < c_{3}\Gamma_{0}^2\Big({{r-r_+}\over 10}\Big)\,\Gamma_{0}^2\Big({{m-m_+}\over 10}\Big)\,{\varepsilon\over {s-s_+}}\\[0.2cm]
& < c_{4}\Gamma_{0}^2\Big({{r-r_+}\over 10}\Big)\,\Gamma_{0}^2\Big({{m-m_+}\over 10}\Big)\,{\varepsilon\over {s}},
\end{align*}
where $c_{4}>c_{3}$, we also used that $3s_+<s.$ With this constant $c_{4}$ we define
$$\vartheta=c_{4}\Gamma_{0}^2\Big({{r-r_+}\over 10}\Big)\,\Gamma_{0}^2\Big({{m-m_+}\over 10}\Big)\,{\varepsilon\over {s}}$$
and rewrite these inequalities in the form
\begin{equation}\label{c23}
\left\{\begin{array}{ll}
\|u\|_{m_+,r_+,s_+}+\|v\|_{m_+,r_+,s_+} <\vartheta s,\\[0.3cm]
\|{{\partial u}\over {\partial \xi}}\|_{m_+,r_+,s_+}+\|{{\partial v}\over {\partial \xi}}\|_{m_+,r_+,s_+}<\vartheta{s\over{r-r_+}}<\vartheta,\\[0.3cm] \|{{\partial u}\over {\partial \eta}}\|_{m_+,r_+,s_+}+\|{{\partial v}\over {\partial \eta}}\|_{m_+,r_+,s_+}<\vartheta,
\end{array}\right.
\end{equation}
whereby we have made use of $s<{{r-r_+}\over 4}.$

The contribution from the functions $u, v$ on the right hand side of (\ref{c22}) can be estimated using the mean value theorem followed by (\ref{c23}), yielding
\begin{align*}
&\|u(\xi+\alpha, \eta) - u(\xi_1,\eta_1)\|_{m_+,r_+,s_+}\\[0.2cm]
&\leq \left\|{{\partial u}\over {\partial \xi}}\right\|_{m_+,r_+,s_+} (|\eta-\alpha|+\|f_+\|_{m_+,r_+,s_+})+\left\|{{\partial u}\over {\partial \eta}}\right\|_{m_+,r_+,s_+}\|g_+\|_{m_+,r_+,s_+}\\[0.2cm]
&<\vartheta{s\over{r-r_+}}|\eta-\alpha|+\vartheta(\|f_+\|_{m_+,r_+,s_+}+\|g_+\|_{m_+,r_+,s_+}).
\end{align*}
The same final estimate is obtained also for the corresponding contribution $v$.

Recalling that $\|f\|_{m,r,s}+\|g\|_{m,r,s}=\varepsilon$, we can use Cauchy's estimate to bound the derivatives of $f,g$ by ${{2\varepsilon}\over s}$ in $D(r_+,s_+)$, so that, again applying the mean value theorem followed by (\ref{c23}), we obtain
\begin{align*}
\|f(\xi + u, \eta + v) -f(\xi,\eta)\|_{m_+,r_+,s_+}&<{{2\varepsilon}\over s}(\|u\|_{m_+,r_+,s_+}+\|v\|_{m_+,r_+,s_+})\\[0.2cm]
&<2\vartheta\varepsilon
\end{align*}
with the same final estimate for the corresponding contribution from $g$.

The troublesome mean value  $[g](\eta)$ will be approximated by the linear function
$$h(\eta)=[g](\alpha)+[g]_{\eta}(\alpha)(\eta-\alpha),$$
which we will estimate later using the intersection property. From (\ref{c22}) and the previous estimates we now have
\begin{align*}
\|f_+\|_{m_+,r_+,s_+}+\|g_+ -h\|_{m_+,r_+,s_+}
&<2\vartheta \Big(\|f_+\|_{m_+,r_+,s_+}+\|g_+ -h\|_{m_+,r_+,s_+}\Big)\\[0.2cm]
&+2\vartheta |h|_{s_+} +2\vartheta{s\over{r-r_+}}|\eta-\alpha|+4\vartheta\varepsilon+\big|[g]-h\big|_{s_+}.
\end{align*}

Now if one chooses $\varepsilon$ sufficiently small such that
$$\vartheta=c_{4}\Gamma_{0}^2\Big({{r-r_+}\over 10}\Big)\,\Gamma_{0}^2\Big({{m-m_+}\over 10}\Big)\,{\varepsilon\over {s}}<{1\over 4},$$
we can eliminate $\|f_+\|_{m_+,r_+,s_+}+\|g_+ -h\|_{m_+,r_+,s_+}$ from the right hand side,
and recalling that $|\eta-\alpha|<s_+<s$, we can express this in the form
\begin{equation*}
\|f_+\|_{m_+,r_+,s_+}+\|g_+ -h\|_{m_+,r_+,s_+}<c_5\Bigg\{\vartheta{s\over{r-r_+}}s+\vartheta |h|_{s_+}+\vartheta\varepsilon+\big|[g]-h\big|_{s_+}\Bigg\}.
\end{equation*}

A preliminary estimate of $h(\eta)$ and $\big|[g]-h\big|_{s_+}$ is obtained by observing that for $|\eta-\alpha|<s$ we have $\Big|[g]\Big|_{s}<\varepsilon$ and therefore by Cauchy's estimate $\Big|[g]_\eta(\alpha)\Big|_{s_+}<{\varepsilon\over{s-s_+}},$ while for $|\eta-\alpha|<s_+$ also $\Big|[g]_{\eta \eta}\Big|_{s_+}<{{2\varepsilon}\over{(s-s_+)^2}}.$ Consequently,  for $|\eta-\alpha|<s_+$, we have
$$|h|_{s_+}<\varepsilon+{\varepsilon\over{s-s_+}}s_+<2\varepsilon$$
and
$$\big|[g]-h\big|_{s_+}<{{s_+^2}\over 2}\Big|[g]_{\eta \eta}\Big|_{s_+}<\Big({{s_+}\over{s-s_+}}\Big)^2 \varepsilon<3\Big({{s_+}\over s}\Big)^2 \varepsilon,$$
where we have used that $3s_+<s.$ Combining these with the previous
estimates, we now have
\begin{align*}
\|f_+\|_{m_+,r_+,s_+}+\|g_+ -h\|_{m_+,r_+,s_+}&<c_5\Bigg\{\vartheta{s\over{r-r_+}}s+2\vartheta \varepsilon+\vartheta\varepsilon+3\Big({{s_+}\over s}\Big)^2 \varepsilon\Bigg\}\\[0.2cm]
&=c_5\Bigg\{{\vartheta\over{r-r_+}}\big(s^2+3(r-r_+)\varepsilon\big)+3\Big({{s_+}\over s}\Big)^2 \varepsilon\Bigg\}\\[0.2cm]
&\leq c_6\Bigg\{{\vartheta\over{r-r_+}}\big(s^2+\varepsilon\big)+3\Big({{s_+}\over s}\Big)^2 \varepsilon\Bigg\},
\end{align*}
this becomes
\begin{align}\label{c17}
& \|f_+\|_{m_+,r_+,s_+}+\|g_+-h\|_{m_+,r_+,s_+}\nonumber\\[0.2cm]
&< c_6\Bigg\{{{c_{4}}\over{r-r_+}} \Gamma_{0}^2\Big({{r-r_+}\over 10}\Big)\,\Gamma_{0}^2\Big({{m-m_+}\over 10}\Big)\Big(s\varepsilon+{{\varepsilon^2}\over s}\Big)+\Big({{s_+}\over s}\Big)^2 \varepsilon\Bigg\}=Q,
\end{align}
where the right hand side of this inequality also deserves to define $Q$.

The preliminary estimate of $2\varepsilon$ for $|h|_+$, however, is insufficient for decreasing the error term, and to obtain a better estimate we use the
intersection property of $\mathfrak{M}$, or that of $\mathfrak{N}=\mathfrak{{U}}^{-1} \mathfrak{M}\, \mathfrak{U}$\, also has the intersection property.\ Accordingly, each curve  $\eta=\text{constant}$, in particular, has to intersect its image curve under $\mathfrak{N}$ and at such a point of intersection we have $\eta_1=\eta$ or $g_+=0$, so that for each real $\eta$ in $|\eta-\alpha|<s_+$ there exists a real $\xi=\xi_0(\eta)$ such that $g_+(\xi_0(\eta),\eta)=0.$ Applying (\ref{c17}) at such points $(\xi_0(\eta),\eta)$, we find that
$$|h(\eta)|<Q\ \ \ (\alpha-s_+<\eta<\alpha+s_+).$$
Consequently, setting $\eta=\alpha$ in the definition of $h$, we get $\big|[g](\alpha)\big|<Q,$ and
letting  $\eta$ approach $\alpha+s_+$ in the same definition, we obtain
$$\big|[g](\alpha)+[g]_{\eta}(\alpha)s_+\big|\leq Q,$$
so that
$$\big|[g]_{\eta}(\alpha)s_+\big|< 2Q.$$
From this we conclude that for complex $\eta$ in the disk $|\eta-\alpha|<s_+$ we have
$$\big|h(\eta)\big|\leq \big|[g](\alpha)\big|+\big|[g]_{\eta}(\alpha)\big|\big|\eta-\alpha\big|<3Q,$$
which in view of (\ref{c17}) gives
\begin{align*}
& \|f_+\|_{m_+,r_+,s_+}+\|g_+\|_{m_+,r_+,s_+}\nonumber\\[0.2cm]
&<4Q= c_{7} \Bigg\{{{1}\over{r-r_+}} \Gamma_{0}^2\Big({{r-r_+}\over 10}\Big)\,\Gamma_{0}^2\Big({{m-m_+}\over 10}\Big)\Big(s\varepsilon+{{\varepsilon^2}\over s}\Big)+\Big({{s_+}\over s}\Big)^2 \varepsilon\Bigg\}.
\end{align*}

\subsection{The Iteration Lemma}

The above discussions lead to the following lemma.

\begin{lemma}\label{lem3.1}
Consider a map
\begin{equation*}
\mathfrak{M}:\quad \begin{array}{ll}
\left\{\begin{array}{ll}
x_1=x+y+f(x,y),\\[0.2cm]
y_1=y+g(x,y),\\[0.1cm]
 \end{array}\right.\
\end{array}
\end{equation*}
where $f$ and $g$ are real analytic in the domain $D(r,s)$ and almost periodic in $x$ with the frequency $\omega=(\cdots,\omega_\lambda,\cdots)$. Assume $\omega=(\cdots,\omega_\lambda,\cdots),\alpha$ satisfy (\ref{b9}). Let
$$\varepsilon=\|f\|_{m,r,s}+\|g\|_{m,r,s}.$$
Then there is a constant $c_4>0$ depending on $\gamma,\Delta$ only such that for
\begin{equation}\label{c24}
0<m_+<m<1,0<r_+<r<1,0<3s_+<s<{{r-r_+}\over 4},
\end{equation}
if
\begin{equation}\label{c12}
\vartheta=c_{4} \Gamma_{0}^2\Big({{r-r_+}\over 10}\Big)\,\Gamma_{0}^2\Big({{m-m_+}\over 10}\Big)\,{\varepsilon\over {s}}<{1\over 4},
\end{equation}
there is a transformation
\begin{equation*}
\begin{array}{ll}
\mathfrak{U}:\ \ \left\{\begin{array}{ll}
x = \xi+u(\xi,\eta),\\[0.2cm]
y = \eta+v(\xi,\eta),
\end{array}\right.
\end{array}
\end{equation*}
which is defined in the domain $D(r_+,s_+)$, $u,v$ are real analytic and  almost periodic in $\xi$. Under this transformation, the map $\mathfrak{M}$ is into the form
\begin{equation*}
\begin{array}{ll}
\mathfrak{{U}}^{-1} \mathfrak{M}\, \mathfrak{U}\ :\ \left\{\begin{array}{ll}
\xi_1=\xi+\eta+f_{+}(\xi,\eta),\\[0.2cm]
\eta_1=\eta+g_{+}(\xi,\eta),
 \end{array}\right.
\end{array}
\end{equation*}
where the functions $f_{+}$ and $g_{+}$ are real analytic in a smaller domain $D(r_{+},s_{+})$ and almost periodic in $\xi$ with the frequency $\omega=(\cdots,\omega_\lambda,\cdots)$. Moreover, the following estimates hold :
\begin{align}
&\|\mathfrak{U}-id\|_{m_+,r_+,s_+}<\vartheta s,\ \ \  \|\partial\mathfrak{U}-E\|_{m_+,r_+,s_+}<\vartheta,\label{c103}\\[0.2cm]
& \|f_+\|_{m_+,r_+,s_+}+\|g_+\|_{m_+,r_+,s_+}\nonumber\\[0.2cm]
&< c_{7} \Bigg\{{{1}\over{r-r_+}} \Gamma_{0}^2\Big({{r-r_+}\over 10}\Big)\,\Gamma_{0}^2\Big({{m-m_+}\over 10}\Big)\Big(s\varepsilon+{{\varepsilon^2}\over s}\Big)+\Big({{s_+}\over s}\Big)^2 \varepsilon\Bigg\}.\label{c14}
\end{align}
\end{lemma}

\section{Iteration and Proof of Theorem \ref{thm2.11}}

Lemma \ref{lem3.1} is usually called the iteration lemma in the KAM proof. For the proof of Theorem \ref{thm2.11}, one can use this lemma infinite times to construct a sequence of transformations $\mathfrak{U}$. To this
end we make a sequence of successive applications of the lemma, starting
with the given mapping (\ref{M}), now denote by $\mathfrak{M}=\mathfrak{M}_0$ and restricted to
the domain
$$D_0\ :\ |\Im\,x|<r_0,\ \ \ |y-\alpha|<s_0.$$
By assumption we have there
$$\|f\|_{m_0,r_0,s_0}+\|g\|_{m_0,r_0,s_0}<\varepsilon_0$$
and therefore have to verify first the inequalities (\ref{c24}), (\ref{c12}). Firstly, we choose
\begin{equation}\label{c25}
\begin{array}{ll}
m_1={3\over 4}m_0,\quad r_1={3\over 4}r_0,\quad s_0=\varepsilon_0^{2\over 3}, \quad s_1=\varepsilon_1^{2\over 3},\\[0.2cm]
\varepsilon_1=r_0^{-1}c_8\Gamma_{0}^{-2}(c_8)\varepsilon_0^{4\over 3}, \quad e_0=r_0^{-3}c_8^{12}\Gamma_{0}^{-6}(c_8)\varepsilon_0
\end{array}
\end{equation}
with $c_8>2$ a suitably chosen constant. Then we have
$$0<m_1={3\over 4}m_0<m_0,\quad 0<r_1={3\over 4}r_0<r_0.$$
The inequality $3s_1<s_0$ follows from
$$\Big({{s_{1}}\over {s_0}}\Big)^{3\over 2}={{\varepsilon_1}\over {\varepsilon_0}}=c_8^{-3}e_0^{1\over 3}<c_8^{-3}<{1\over 8}$$
and since $ s_0=\varepsilon_0^{2\over 3}$ goes to $0$ faster than exponentially when $\varepsilon_0$ is small, while
$r_0-r_1=r_0 2^{-2}$ decays only exponentially, the whole third inequality in (\ref{c24}) will hold
for $\varepsilon_0$ sufficiently small.
And for $r=r_0,r_+=r_1,m=r_0,m_+=m_1,\varepsilon=\varepsilon_0,s=s_0$ in (\ref{c12}) we have
$$\vartheta_0=c_{4} \Gamma_{0}^2\Big({{r_0}\over 40}\Big)\,\Gamma_{0}^2\Big({{m_0}\over 40}\Big)\,{{\varepsilon_0}\over {s_0}}=c_{4} \Gamma_{0}^2\Big({{r_0}\over 40}\Big)\,\Gamma_{0}^2\Big({{m_0}\over 40}\Big)\,{\varepsilon_0}^{1\over 3},$$
this can be made less than ${1\over 4}$ by choosing  $\varepsilon_0$ small.

Hence ${\varepsilon}_0$ can be chosen sufficiently small such that the inequalities (\ref{c24}), (\ref{c12}) with $n=0$ hold. Transforming the mapping $\mathfrak{M}_0$ by
the coordinate transformation $\mathfrak{U}=\mathfrak{U}_0$ provided by Lemma  \ref{lem3.1}, we
obtain a mapping  $\mathfrak{M}_1=\mathfrak{{U}}^{-1} \mathfrak{M}\, \mathfrak{U}$ defined in the domain
$$D_1\ :\ |\Im\,x|<r_1,\ \ \ |y-\alpha|<s_1,$$
where $r_1,s_1$ correspond to the parameter $r_+,s_+$ of Lemma  \ref{lem3.1}. Applying Lemma  \ref{lem3.1} to the new mapping $\mathfrak{M}_1$, we obtain another coordinate transformation $\mathfrak{U}_1$ and a transformed mapping  $\mathfrak{M}_2={\mathfrak{U}_1}^{-1} \mathfrak{M}_1\, \mathfrak{U}_1$, and
proceeding in this way we are led to a sequence of mappings
\begin{equation}\label{c15}
\mathfrak{M}_{n+1}={\mathfrak{U}_n}^{-1} \mathfrak{M}_n\, \mathfrak{U}_n\ \ \ (n=0,1,\cdots)
\end{equation}
whose domains $D_{n+1}$ are defined like $D_1$ with $r_{n+1},s_{n+1}$ replacing $r,s$. We have to verify, of course, that this sequence of transformations are well defined, and that $\mathfrak{M}_n$ approximates the twist mapping with
increasing precision. For this we fix the parameters $m_n,r_n,s_n,\varepsilon_n\,(n=0,1,\cdots)$ by setting
\begin{equation}\label{c16}
\begin{array}{ll}
m_n={{m_0}\over 2}\Big(1+{1\over{2^n}}\Big),\ \ \ \ r_n={{r_0}\over 2}\Big(1+{1\over{2^n}}\Big),\\[0.2cm]
s_{n}=\varepsilon_n^{2\over 3},\ \ \ \  \varepsilon_{n+1}=r_0^{-1}c_8^{n+1}\Gamma_{0}^{-2}(c_8^{n+1})\varepsilon_n^{4\over 3},
\end{array}
\end{equation}
where $c_8$ is given in (\ref{c25}). Thus $m_n$ is a decreasing sequence converging to the positive value ${{m_0}\over 2}$, $r_n$ is also a decreasing sequence converging to the positive value ${{r_0}\over 2}$, and all functions to be considered will be analytic in $\xi$ for $|\Im\,\xi|<{{r_0}\over 2}$. The sequence $\varepsilon_n$ converges to $0$ provided $\varepsilon_0$ chosen sufficiently small. Indeed, the sequence $e_n=r_0^{-3}c_8^{3(n+4)}\Gamma_{0}^{-6}(c_8^{n+1})\varepsilon_n$  satisfies
$$e_{n+1}\leq e_n^{4\over 3}$$
and therefore converges to zero if we take $0<e_0<1$, or $0<\varepsilon_0<r_0^{3}c_8^{-12}\Gamma_{0}^{6}(c_8)$. Then $s_n$ is also a decreasing sequence converging to $0$.

To show that the mapping $\mathfrak{M}_n$ is well defined in $D_n$ and satisfies the appropriate estimate, we proceed by induction. Assume that  $\mathfrak{M}_n$ is defined in $D_n$ and satisfies the estimate
$$\varepsilon=\|f\|_{m_n,r_n,s_n}+\|g\|_{m_n,r_n,s_n}<\varepsilon_n,$$
we will verify the corresponding statement for $\mathfrak{M}_{n+1}$. For this we apply Lemma \ref{lem3.1} with $m=m_n, r=r_n, s=s_n, m_{+}=m_{n+1}, r_{+}=r_{n+1}, s_{+}=s_{n+1},$  and therefore have to verify first the inequalities (\ref{c24}), (\ref{c12}). The inequality $3s_{n+1}<s_n$ follows from
\begin{equation}\label{c26}
\Big({{s_{n+1}}\over {s_n}}\Big)^{3\over 2}={{\varepsilon_{n+1}}\over {\varepsilon_n}}=c_8^{-3}e_n^{1\over 3}<c_8^{-3}<{1\over 8}
\end{equation}
and since $\varepsilon_n$ goes to $0$ faster than exponentially while
$r-r_+=r_n-r_{n+1}=r_0 2^{-n-2}$ decays only exponentially, the whole third inequality in (\ref{c24}) will hold
for $\varepsilon_0$ sufficiently small.
And for $\vartheta=\vartheta_n$ in (\ref{c12}) we have
$$\vartheta_n=c_{4} \Gamma_{0}^2\Big({{r_n-r_{n+1}}\over 10}\Big)\,\Gamma_{0}^2\Big({{m_n-m_{n+1}}\over 10}\Big)\,\varepsilon_n^{1\over 3 },$$
this can be made less than ${1\over 4}$ by choosing  $\varepsilon_0$ small. Thus there exists a positive constant $\varepsilon^*=\varepsilon^*(r_0,s_0,m_0,\gamma,\Delta)$ such that for $\varepsilon_0<\varepsilon^*$  the inequalities (\ref{c24}),(\ref{c12}) hold and Lemma \ref{lem3.1} is applicable. From Lemma \ref{lem3.1} we now obtain the transformation ${\mathfrak{U}_n}$ taking $D_{n+1}$ into $D_n$ and
the transformed mapping $\mathfrak{M}_{n+1}={\mathfrak{U}_n}^{-1} \mathfrak{M}_n\, \mathfrak{U}_n$ defined in $D_{n+1}$. Moreover, by (\ref{c14}) we have the estimate
\begin{align*}
&\|f_+\|_{m_{n+1},r_{n+1},s_{n+1}}+\|g_+\|_{m_{n+1},r_{n+1},s_{n+1}}\nonumber\\[0.2cm]
&<  c_{7} \Bigg\{{{1}\over{r_n-r_{n+1}}} \Gamma_{0}^2\Big({{r_n-r_{n+1}}\over 10}\Big)\,\Gamma_{0}^2\Big({{m_n-m_{n+1}}\over 10}\Big)\,\Big(\varepsilon_n^{5\over 3}+\varepsilon_n^{4\over 3}\Big)+\Big({{\varepsilon_{n+1}}\over {\varepsilon_n}}\Big)^{4\over 3} \varepsilon_n\Bigg\}\nonumber\\[0.2cm]
&=  c_{7} \Bigg\{{{1}\over{r_n-r_{n+1}}} \Gamma_{0}^2\Big({{r_n-r_{n+1}}\over 10}\Big)\,\Gamma_{0}^2\Big({{m_n-m_{n+1}}\over 10}\Big)\,\varepsilon_n^{4\over 3}\Big(\varepsilon_n^{1\over 3}+1\Big)+\Big({{\varepsilon_{n+1}}\over {\varepsilon_n}}\Big)^{1\over 3} \varepsilon_{n+1}\Bigg\},
\end{align*}
which in view of (\ref{c16}), (\ref{c26}) gives
\begin{align*}
&\|f_+\|_{m_{n+1},r_{n+1},s_{n+1}}+\|g_+\|_{m_{n+1},r_{n+1},s_{n+1}}\nonumber\\[0.1cm]
&\leq c_{7} \Bigg\{2^{n+2} c_8^{-n-1}\Gamma_{0}^{2}(c_8^{n+1})\Gamma_{0}^2\Big({{r_0 2^{-n-2}}\over 10}\Big)\,\Gamma_{0}^2\Big({{m_0 2^{-n-2}}\over 10}\Big)\,\Big(\varepsilon_n^{1\over 3}+1\Big)+\Big({{\varepsilon_{n+1}}\over {\varepsilon_n}}\Big)^{1\over 3} \Bigg\}\,\varepsilon_{n+1}\nonumber\\[0.1cm]
&<  c_{7} \Bigg\{2(2 c_8^{-1})^{n+1}\Gamma_{0}^{2}(c_8^{n+1})\Gamma_{0}^2\Big({{r_0 2^{-n-2}}\over 10}\Big)\,\Gamma_{0}^2\Big({{m_0 2^{-n-2}}\over 10}\Big)\,\Big(\varepsilon_n^{1\over 3}+1\Big)+c_8^{-1}\Bigg\}\, \varepsilon_{n+1},
\end{align*}
since the sequence $\varepsilon_n$ is bounded, and in view of condition (\ref{b15}) the coefficient of $\varepsilon_{n+1}$ in the last
inequality can be made less than $1$ by taking $c_8$ large. With such a choice
we finally have
\begin{equation*}
\|f_+\|_{m_{n+1},r_{n+1},s_{n+1}}+\|g_+\|_{m_{n+1},r_{n+1},s_{n+1}}<\varepsilon_{n+1}
\end{equation*}
and the induction is complete.

Since ${\mathfrak{U}_k}$ maps the domain $D_{k+1}$ into $D_{k}\,(k=0,1,2,\cdots)$, the
transformation $\mathfrak{V}_n={\mathfrak{U}_0}\circ{\mathfrak{U}_1}\circ\cdots{\mathfrak{U}_n}$ is well defined in $D_{n+1}$ and is seen to take $\mathfrak{M}_{0}$ into
$$\mathfrak{M}_{n+1}={\mathfrak{V}_n}^{-1}\mathfrak{M}_{0}{\mathfrak{V}_n}.$$
Moreover, if we express ${\mathfrak{V}_n}$ in the form
$$x = \xi+p_n(\xi,\eta),y = \eta+q_n(\xi,\eta),$$
then $p_n,q_n$ are analytic functions in the domain $D_{n+1}$ which, in view of $s_n\rightarrow 0, r_n\rightarrow {{r_0}\over 2},$ shrinks down to the domain $D_\infty$ :
$$|\Im\,\xi|<{{r_0}\over 2},\ \ |\eta-\alpha|=0.$$
We will show that as $n\rightarrow \infty$ the sequences $p_n(\xi,\alpha),q_n(\xi,\alpha)$ converge to
analytic functions of $\xi$ for $|\Im\,\xi|<{{r_0}\over 2}$. Indeed, from $\mathfrak{V}_n=\mathfrak{V}_{n-1}{\mathfrak{U}_n}$ it
follows that
\begin{equation}\label{c102}
\begin{array}{ll}
\begin{array}{ll}
p_n = u_n+p_{n-1}(\xi+u_n,\eta+v_n),\\[0.2cm]
q_n = v_n+q_{n-1}(\xi+u_n,\eta+v_n),
\end{array}
\end{array}
\end{equation}
where $u_n,v_n$ correspond to $\mathfrak{U}_n$ as $u,v$ to $\mathfrak{U}$, and consequently
\begin{equation*}
\begin{array}{ll}
\begin{array}{ll}
p_n = u_n+u_{n-1}+\cdots+u_0,\\[0.2cm]
q_n = v_n+v_{n-1}+\cdots+v_0.
\end{array}
\end{array}
\end{equation*}
In the last two sums we have \mbox{supp}ressed the different arguments, but this
is irrelevant for the proof of convergence since we can estimate $\|u_n\|_{m_{n+1},r_{n+1},s_{n+1}},\\ \|v_n\|_{m_{n+1},r_{n+1},s_{n+1}}$ by
their supremum over $D_{n+1}$. Namely, from (\ref{c103}) it follows that
$$\|u_n\|_{m_{n+1},r_{n+1},s_{n+1}}+\|v_n\|_{m_{n+1},r_{n+1},s_{n+1}}<{1\over 4}s_n,$$
Also, from $3s_{n+1}<s_n$ we have
\begin{equation}\label{c100}
\|p_n\|_{m_{n+1},r_{n+1},s_{n+1}}+\|q_n\|_{m_{n+1},r_{n+1},s_{n+1}}<{1\over 4}\sum_{\nu=0}^{\nu=n} s_\nu<{1\over 4}\sum_{\nu=0}^{\infty}{1\over{3^\nu}} s_0< s_0.
\end{equation}

Estimate the derivatives of $p_n,q_n$ on $D_{n+1}$. Let
$$P_n^1=\max\Bigg\{\|{{\partial {p_n}}\over {\partial \xi}}\|_{m_{n+1},r_{n+1},s_{n+1}}, \|{{\partial {q_n}}\over {\partial \xi}}\|_{m_{n+1},r_{n+1},s_{n+1}}\Bigg\},$$
$$P_n^2=\max\Bigg\{\|{{\partial {p_n}}\over {\partial \eta}}\|_{m_{n+1},r_{n+1},s_{n+1}}, \|{{\partial {q_n}}\over {\partial \eta}}\|_{m_{n+1},r_{n+1},s_{n+1}}\Bigg\}.$$
Show by induction that for all $n\geq 0$
\begin{equation}\label{c101}
P_n^1+P_n^2<{\prod \limits _{\ell=0}^n }(1+2\vartheta_\ell)-1.
\end{equation}

For $n=0$, the inequality (\ref{c101}) is fulfilled, since
$$\|{{\partial u}\over {\partial \xi}}\|_{m_{n+1},r_{n+1},s_{n+1}}+\|{{\partial v}\over {\partial \xi}}\|_{m_{n+1},r_{n+1},s_{n+1}}<\vartheta_n,$$
$$\|{{\partial u}\over {\partial \eta}}\|_{m_{n+1},r_{n+1},s_{n+1}}+\|{{\partial v}\over {\partial \eta}}\|_{m_{n+1},r_{n+1},s_{n+1}}<\vartheta_n$$
on $D_{n+1}$ for all $n$ since which is guaranteed by (\ref{c103}).

Let $n\geq 1$. Since ${\mathfrak{U}_n}(D_{n+1})\subseteq D_{n}$, and by (\ref{c102}) implies
$$\|{{\partial {p_n}}\over {\partial \xi}}\|_{m_{n+1},r_{n+1},s_{n+1}}<\vartheta_n+P_{n-1}^1(1+\vartheta_n)+P_{n-1}^2\vartheta_n,$$
$$\|{{\partial {q_n}}\over {\partial \xi}}\|_{m_{n+1},r_{n+1},s_{n+1}}<\vartheta_n+P_{n-1}^1(1+\vartheta_n)+P_{n-1}^2\vartheta_n,$$
$$\|{{\partial {p_n}}\over {\partial \eta}}\|_{m_{n+1},r_{n+1},s_{n+1}}<\vartheta_n+P_{n-1}^1 \vartheta_n +P_{n-1}^2(1+\vartheta_n),$$
$$\|{{\partial {q_n}}\over {\partial \eta}}\|_{m_{n+1},r_{n+1},s_{n+1}}<\vartheta_n+P_{n-1}^1 \vartheta_n+P_{n-1}^2(1+\vartheta_n).$$

Thus
\begin{align*}
P_n^1+P_n^2 &<(1+2\vartheta_n)(P_{n-1}^1+P_{n-1}^2)+2\vartheta_n\\[0.2cm]
&<(1+2\vartheta_n)\Bigg\{{\prod \limits _{\ell=0}^{n-1} }(1+2\vartheta_{\ell})-1\Bigg\}+2\vartheta_n\\[0.2cm]
&<{\prod \limits _{\ell=0}^n }(1+2\vartheta_\ell)-1.
\end{align*}
We have proved the inequality (\ref{c101}) for all $n$. Note that (\ref{c101}) implies
\begin{equation}\label{c1000}
P_n^1+P_n^2<\Big({3\over 2}\Big)^n-1
\end{equation}
for all $n$.

And by (\ref{c102}), we have
\begin{equation*}
\begin{array}{ll}
\begin{array}{ll}
p_n-p_{n-1} = u_n+p_{n-1}(\xi+u_n,\eta+v_n)-p_{n-1}(\xi,\eta),\\[0.2cm]
q_n-q_{n-1} = v_n+q_{n-1}(\xi+u_n,\eta+v_n)-q_{n-1}(\xi,\eta),
\end{array}
\end{array}
\end{equation*}
then
\begin{align*}
&\|p_n-p_{n-1}\|_{m_{n+1},r_{n+1},s_{n+1}}\\[0.0001cm]
&\leq\|u_n\|_{m_{n+1},r_{n+1},s_{n+1}}+P_{n-1}^1 \|u_n\|_{m_{n+1},r_{n+1},s_{n+1}}+P_{n-1}^2 \|v_n\|_{m_{n+1},r_{n+1},s_{n+1}}\\[0.0001cm]
&\leq(\|u_n\|_{m_{n+1},r_{n+1},s_{n+1}}+\|v_n\|_{m_{n+1},r_{n+1},s_{n+1}})(1+P_{n-1}^1+P_{n-1}^2),\\[0.0001cm]
&\|q_n-q_{n-1}\|_{m_{n+1},r_{n+1},s_{n+1}}\\[0.0001cm]
&\leq\|v_n\|_{m_{n+1},r_{n+1},s_{n+1}}+P_{n-1}^1 \|u_n\|_{m_{n+1},r_{n+1},s_{n+1}}+P_{n-1}^2 \|v_n\|_{m_{n+1},r_{n+1},s_{n+1}}\\[0.0001cm]
&\leq(\|u_n\|_{m_{n+1},r_{n+1},s_{n+1}}+\|v_n\|_{m_{n+1},r_{n+1},s_{n+1}})(1+P_{n-1}^1+P_{n-1}^2).
\end{align*}
By (\ref{c103}), we have
$$\|u_n\|_{m_{n+1},r_{n+1},s_{n+1}}+\|v_n\|_{m_{n+1},r_{n+1},s_{n+1}}<\vartheta_n s_n.$$
And ${\mathfrak{U}_n}(D_{\infty})\subseteq{\mathfrak{U}_n}(D_{n+1})\subseteq D_{n}$, therefore by (\ref{c101}),(\ref{c1000}),
$$\|p_n-p_{n-1}\|_{{m_{0}\over 2},{r_{0}\over 2}}<\vartheta_n s_n(1+P_{n-1}^1+P_{n-1}^2)<{{s_n}\over 4}\Big({3\over 2}\Big)^n,$$
$$\|q_n-q_{n-1}\|_{{m_{0}\over 2},{r_{0}\over 2}}<\vartheta_n s_n(1+P_{n-1}^1+P_{n-1}^2)<{{s_n}\over 4}\Big({3\over 2}\Big)^n.$$

One sees
$$\sum \limits _{n=1}^\infty \Big({3\over 2}\Big)^n s_n<\infty,$$
since
$s_n<{{r_n-r_{n+1}}\over 4}=r_02^{-n-4}$ which is guaranteed by (\ref{c24}), in other words $s_{n}=\varepsilon_n^{2\over 3}$ decay faster than a geometrical progression. Thus this implies uniform convergence of $p_n,q_n$ in $|\Im\,\xi|<{{r_0}\over 2}$ for $\eta=\alpha.$ We denote the limit of $p_n$ by $\widetilde{u}(\xi)$ and of $q_n$ by $\widetilde{v}(\xi)-\alpha$, which then are real analytic functions in  $|\Im\,\xi|<{{r_0}\over 2}$ and almost periodic in $\xi$ with the frequency $\omega=(\cdots,\omega_\lambda,\cdots)$. And we choose $\varepsilon_0= s_0^{3\over 2}$ so small that $ s_0<\bar{\varepsilon}$, therefore from (\ref{c100}), there exists a positive $\tilde{\varepsilon}_0=\tilde{\varepsilon}_0(\bar{\varepsilon},r_0,s_0,m_0,\gamma,\Delta)<\varepsilon^*$ such that for $\varepsilon_0<\tilde{\varepsilon}_0$ we have  $$\|\widetilde{u}\|_{{m_{0}\over 2},{r_{0}\over 2}}+\|\widetilde{v}-\alpha\|_{{m_{0}\over 2},{r_{0}\over 2}}<\bar{\varepsilon}$$
in $|\Im\,\xi|<{{r_0}\over 2}$, as asserted in the theorem, and this $\tilde{\varepsilon}_0$ can be used as the
constant restricting the size of $\|f\|_{m,r,s}+\|g\|_{m,r,s}$. This concludes the proof of the existence of an invariant curve. Hence the proof of Theorem \ref{thm2.11} is completed.\qed

\section{The small twist theorem}

In this section we formulate an useful  small twist theorem which is a variant of the invariant curve theorem (Theorem \ref{thm2.11}) for the almost periodic  mapping  $\mathfrak{M}$.

In many applications,\ one may meet the so called small twist mappings
\begin{equation}\label{f1}
\begin{array}{ll}
\mathfrak{M}_{\delta}:\ \ \left\{\begin{array}{ll}
x_1=x+\delta y+f(x,y,\delta),\\[0.2cm]
y_1=y+g(x,y,\delta),\\[0.1cm]
 \end{array}\right.\  (x,y)\in \mathbb{R} \times [a,b],
\end{array}
\end{equation}
where the functions $f$ and $g$ are real analytic functions, and  almost periodic in $x$ with the frequency $\omega=(\cdots,\omega_\lambda,\cdots)$, $0<\delta< 1$ is a small parameter.

We choose the number $\alpha$ satisfying the inequalities
\begin{equation}\label{f2}
\begin{array}{ll}
\left\{\begin{array}{ll}
a+\gamma\leq \alpha\leq b-\gamma,\\[0.4cm]
\Big|{\langle k,\omega \rangle {{\delta\alpha} \over {2\pi}}-j}\Big|\geq {\gamma \over {\Delta([[k]])\Delta(|k|)}},\ \ \ \  \mbox{for all}\ \  k \in \mathbb{Z}_{0}^{{\mathbb{Z}}}\backslash\{0\},\ \  j \in \mathbb{Z}
 \end{array}\right.
\end{array}
\end{equation}
with some positive constant $\gamma$, where $\Delta$ is some approximation functions.

\begin{theorem}\label{thm6.1}
Suppose that the  almost periodic mapping $\mathfrak{M}_{\delta}$ given by (\ref{f1}) has the intersection property, and  for every $y,\delta$, $f(\cdot,y,\delta),g(\cdot,y,\delta)\in AP_r(\omega)$ with $\omega$ satisfying the nonresonance condition (\ref{b6}), and the corresponding shell functions $F(\theta,y,\delta),G(\theta,y,\delta)$ are real
analytic in the domain $D(r,s)=\big\{(\theta,y)\in \mathbb{C}^\mathbb{Z} \times \mathbb{C}\ :\ |\Im\,\theta|_{\infty}<r, |y-\alpha|<s\big\}$ with $\alpha$ satisfying (\ref{f2}). Then for each positive $\bar{\varepsilon}$ there is a positive $\varepsilon_0=\varepsilon_0(\bar{\varepsilon},r,s,m,\gamma,\Delta)$ such that if $f,g$ satisfy  the following smallness condition
\begin{equation}\label{f3}
\|f(\cdot,\cdot,\delta)\|_{m,r,s}+\|g(\cdot,\cdot,\delta\|_{m,r,s}<\delta{{\varepsilon}}_0,
\end{equation}
then the almost periodic mapping $\mathfrak{M}_{\delta}$ has an invariant curve $\mathbf{\Gamma _{0}}$ with the form
$$
\begin{array}{ll}
\left\{\begin{array}{ll}
x=x'+\varphi(x'),\\[0.2cm]
y=\psi(x'),
 \end{array}\right.\ \ \
\end{array}
$$
where $\varphi, \psi$ are almost periodic  with frequencies $\omega=(\cdots,\omega_\lambda,\cdots)$, and the invariant curve $\mathbf{\Gamma _{0}}$ is of the form $y=\phi(x)$ with $\phi\in AP_{r'}(\omega)$ for some $r'<r$, and $\|\phi-\alpha\|_{m',r'}<\bar{\varepsilon},\ 0<m'<m$. Moreover,\ the  restriction of $\mathfrak{M}_{\delta}$ onto $\mathbf{\Gamma _{0}}$ is of the form
$$\mathfrak{M}_{\delta}|_{\mathbf{\Gamma _{0}}}:\ \ \ \  x_{1}^\prime=x^\prime+\delta\alpha.$$
\end{theorem}

\begin{remark}
If all the conditions of Theorem \ref{thm6.1} hold,  given any $\alpha$ satisfying the inequalities (\ref{f2}), there exists an invariant curve $\mathbf{\Gamma _{0}}$ of $\mathfrak{M}_{\delta}$, which is almost periodic  with the frequency  $\omega=(\cdots,\omega_\lambda,\cdots)$,\ and the restriction of $\mathfrak{M}_{\delta}$ onto $\mathbf{\Gamma _{0}}$ has the form
$$\mathfrak{M}_{\delta}|_{\mathbf{\Gamma _{0}}}:\   x_{1}^\prime=x^\prime +\delta\alpha .$$
\end{remark}

This is the so called small twist theorem.\ It is not a direct consequence of Theorem \ref{thm2.11},\ but one can use the same procedure in the proof of Theorem \ref{thm2.11} to prove it.\ Because there is nothing new in the proof,\ we omit it here.
\begin{remark}
The above conclusion is also true for the following mapping
\begin{equation}\label{f4}
\begin{array}{ll}
\left\{\begin{array}{ll}
x_1=x+\beta+\delta h(y)+f(x,y,\delta),\\[0.2cm]
y_1=y+g(x,y,\delta)
 \end{array}\right.
\end{array}
\end{equation}
with $h{'}(y)\neq 0$, if we change the conditions (\ref{f3}) into
$$
\begin{array}{ll}
M\Big(\|f(\cdot,\cdot,\delta)\|_{m,r,s}+\|g(\cdot,\cdot,\delta)\|_{m,r,s}\Big)<{\varepsilon}_0
\end{array}
$$
with $M=\max\Big\{|h{'}|_{s}\, , 1 \Big\}$.
\end{remark}

\section{Application}

In this section we will apply the above result to superlinear Duffing's equation
\begin{equation}\label{g1}
\ddot{x}+x^3=f(t),
\end{equation}
where $f(t)$ is a real analytic almost periodic function with the frequency $\omega=(\cdots,\omega_\lambda,\cdots)$ and admit a rapidly converging Fourier series expansion.

\subsection{Action and angle variables}

Dropping the time-dependent term equation (\ref{g1}) becomes
$$\ddot{x}+x^3=0.$$
Introduce a new variable $y=\dot{x}$ yields
\begin{equation}\label{g2}
\left\{\begin{array}{ll}
\dot{x}=y, \\[0.1cm]
\dot{y}=-x^3,
\end{array}\right.
\end{equation}
which is time-independent Hamiltonian system on $\mathbb{R}^2$:
\begin{equation}\label{g3}
\left\{\begin{array}{ll}
\dot{x}={\partial\over{\partial y}}h(x,y), \\[0.2cm]
\dot{y}=-{\partial\over{\partial x}}h(x,y) \\
\end{array}\right.
\end{equation}
with $h(x,y)={1\over 2}y^2+{1\over 4}x^4.$ Clearly, $h$ is positive on $\mathbb{R}^2$ except at the only equilibrium point $(x,y)=(0,0)$ of (\ref{g2}) where $h=0$. All the solutions of (\ref{g3}) are periodic with period tending to zero as $h=E$ tends to infinity.

Suppose $(C(t),S(t))$ is the solution of (\ref{g2}) satisfying the initial condition $(C(0),S(0))=(1,0)$. Let $T_*>0$ be its minimal period.
From (\ref{g2}), these analytic functions satisfy\\
(i) $C(t+T_*)=C(t), S(t+T_*)=S(t)\ \text{and}\ C(0)=1, S(0)=0.$\\[0.2cm]
(ii) $\dot{C}(t)=S(t), \dot{S}(t)=-C^3(t).$\\[0.2cm]
(iii) $2S^2(t)+C^4(t)=1.$\\[0.2cm]
(iiii) $C(-t)=C(t), S(-t)=-S(t).$

The action and angle variables are now defined by the map $\Psi :\mathbb{R}^+\times \mathbb{S}^1\rightarrow \mathbb{R}^2\setminus {0}$, where $(x,y)=\Psi(\theta,\rho)$ with $\rho>0$ and $\theta$ $(\text{mod}\ 1)$ is given by the formula
\begin{equation}
\Psi :\begin{array}{ll}
x=c^{1\over 3}{\rho^{1\over 3}}C( {\theta{T_*}}), \\[0.2cm]
y=c^{2\over 3}{\rho^{2\over 3}}S( {\theta{T_*}}), \\[0.2cm]
\end{array}
\end{equation}
with $c={3\over {T_*}}$.

We claim that $\Psi$ is a symplectic diffeomorphism from $\mathbb{R}^+\times \mathbb{S}^1$ onto $\mathbb{R}^2\setminus {0}$. Indeed, for the Jacobian $\tilde{\Delta}$ of $\Psi$ one finds by (ii) and (iii) $|\tilde{\Delta}|=1$, so that $\Psi$ is measure preserving. Moreover since $(C,S)$ is a solution of a differential equation having $T_*$ as minimal period one concludes that $\Psi$ is one to one and onto, which proves the claim.

In the new coordinates  the Hamiltonian associated to (\ref{g2}) becomes
$$h\circ \Psi(\rho,\theta)=d \cdot \rho^{4\over 3}=h_0(\rho),$$
$d={{c^{4\over 3}}\over 4}$, which is independent of the angle variable $\theta$ so that system (\ref{g2}) becomes very simple :
\begin{equation*}
\left\{\begin{array}{ll}
\dot{\theta}={\partial\over{\partial \rho}}h_0={4\over 3}d\rho^{1\over 3}, \\[0.2cm]
\dot{\rho}=-{\partial\over{\partial \theta}}h_0=0.
\end{array}\right.
\end{equation*}

The full equation (\ref{g1}) has the Hamiltonian function :
\begin{equation*}
h(x,y,t)={1\over 2}y^2+{1\over 4}x^4-x f(t).
\end{equation*}
Under the symplectic transformation $\Psi$ it is transformed into :
\begin{equation*}
h_1(\rho,\theta,t)=h(\Psi(\rho,\theta),t)=d \cdot \rho^{4\over 3}-c^{1\over 3}{\rho^{1\over 3}}C( {\theta{T_*}})f(t).
\end{equation*}
The Hamiltonian system (\ref{g1}) now becomes as follows :
\begin{equation*}
\left\{\begin{array}{ll}
\dot{\theta}={\partial\over{\partial \rho}}h_1={4\over 3}d\rho^{1\over 3}-{1\over 3}c^{1\over 3}{\rho^{-{2\over 3}}}C( {\theta{T_*}})f(t), \\[0.2cm]
\dot{\rho}=-{\partial\over{\partial \theta}}h_1={{c^{1\over 3}} T_*}{\rho^{1\over 3}}S( {\theta{T_*}})f(t).
\end{array}\right.
\end{equation*}

As one did in the periodic case, now we change the role of the variable $t$ and $\theta$, and yields that
\begin{equation}\label{g5}
\left\{\begin{array}{ll}
{dt\over d\theta} =\big[{4\over 3}d\rho^{1\over 3}-{1\over 3}c^{1\over 3}{\rho^{-{2\over 3}}}C( {\theta{T_*}})f(t)\big]^{-1}, \\[0.2cm]
{d\rho\over d\theta} ={{c^{1\over 3}} T_*}{\rho^{1\over 3}}S( {\theta{T_*}})f(t)\big[{4\over 3}d\rho^{1\over 3}-{1\over 3}c^{1\over 3}{\rho^{-{2\over 3}}}C( {\theta{T_*}})f(t)\big]^{-1},
\end{array}\right.
\end{equation}
this system is $1$-periodic in the new time variable $\theta$. Let $\rho_{*}$ be a positive number such that
$${4\over 3}d\rho_*^{1\over 3}-{1\over 3}c^{1\over 3}{\rho_*^{-{2\over 3}}}\big|C\big|\big|f\big|>0.$$
System (\ref{g5}) is well defined for $\rho\geq \rho_{*}$. Let $\big(t(\theta),r(\theta)\big)$ be a solution of (\ref{g5}) defined in a certain interval $I=[\theta_{0},\theta_{1}]$ and such that $\rho(\theta)> \rho_{*}$ for all $\theta$ in $I$. The derivative ${dt\over d\theta}$ is positive and the function $t$ is a diffeomorphism from $I$ onto $J=[t_{0},t_{1}]$, where $t(\theta_{0})=t_{0}$ and $t(\theta_{1})=t_{1}$. The inverse function will be denoted by $\theta=\theta(t)$. It maps $J$ onto $I$.

\subsection{The expression of the Poincar\'{e} map of (\ref{g5})}

Our next goal is to obtain asymptotic expansions for $t_{1}$ and $\rho_{1}$. For $\rho_0$ large enough, the second equation of (\ref{g5}) can be rewritten as
$${d\rho\over d\theta}={{3 T_* c^{1\over 3}}\over {4 d}}S( {\theta{T_*}})f(t)+O(\rho^{-1}).$$
An integration of this equation leads to
\begin{equation}\label{g6}
\begin{array}{ll}
\rho(\theta)=\rho_0+O(1),& \theta \in [0,1].
\end{array}
\end{equation}
Therefore,
$$\rho(\theta)^{-{1\over 3}}={{\rho_0}^{-{1\over 3}}}(1+{\rho_0}^{-1}O(1))^{-{1\over 3}}.$$
Expanding $(1+{\rho_0}^{-1}O(1))^{-{1\over 3 }}$ yields
\begin{equation}\label{g7}
\begin{array}{ll}
{\rho(\theta)}^{-{1\over 3}}={\rho_0}^{-{1\over 3}}+O({\rho_0}^{-{4\over 3}}),& \theta \in [0,1]
\end{array}
\end{equation}
for $\rho_0$ large enough,\ by (\ref{g7}) and the first equality of (\ref{g5}),\ we get
$${dt\over d\theta}=O({\rho_0}^{-{1\over 3}}),\ \theta \in [0,1],$$
which implies that
\begin{equation}\label{g8}
\begin{array}{ll}
t(\theta)=t_0+O({\rho_0}^{-{1\over 3}}),& \theta \in [0,1]
\end{array}
\end{equation}
for $\rho_0$ large enough.\ Substituting (\ref{g6}), \ (\ref{g7}),\ (\ref{g8}) into the second equality of (\ref{g5}),\ we have,\ for $\theta \in [0,1]$,
\begin{eqnarray}\label{g9}
{d\rho\over d\theta}={{3 T_* c^{1\over 3}}\over {4 d}}S( {\theta{T_*}})f(t_0+\rho_0^{-{1\over3}})+O(\rho_0^{-1})
\end{eqnarray}
for $\rho_0$ large enough.\ An integration of (\ref{g9}) over $\theta \in [0,1]$ yields
$$
{\rho_1} = \rho_0+O(\rho_0^{-{1\over3}})
$$
for $\rho_0$ large enough, where $\rho_1=\rho(1)$.\ Substituting (\ref{g7}),\ (\ref{g8}) into the first equality of (\ref{g5}),\ we have,\ for $\theta \in [0,1]$
\begin{eqnarray}\label{g11}
{dt\over d\theta}={{3}\over {4d}}{\rho_0}^{-{1\over 3}}+O({\rho_0}^{-{4\over 3}})
\end{eqnarray}
for $r_0$ large enough.\ An integration (\ref{g11}) over $\theta \in [0,1]$ yields
$$
t_1=t_0+{{3}\over {4d}}{\rho_0}^{-{1\over 3}}+O({\rho_0}^{-{4\over 3}}),
$$
where $t_1=t(1)$.

Then the Poincar\'{e} map $P$ of (\ref{g5}) has the expansion
\begin{equation*}
P:
\left\{\begin{array}{ll}
t_1=t_0+{{3}\over {4d}}{\rho_0}^{-{1\over 3}}+O({\rho_0}^{-{4\over 3}}), \\[0.2cm]
{\rho_1} = \rho_0+O(\rho_0^{-{1\over3}}).
\end{array} \right.
\end{equation*}

\subsection{Intersection property of the Poincar\'{e} map $P$ of (\ref{g5})}

Since the Poincar\'{e} map $P$ of (\ref{g5}) ia an exact symplectic map, then $P$ possesses the intersection property which is guaranteed by Lemma \ref{lem7.10}.

To apply the small twist theorem obtained in Section 6, we need a fixed annulus, for this reason, introduce a new variable $\mu$ and a small parameter $\delta >0$ by
\begin{equation}\label{g13}
\delta \mu={{3}\over {4d}}\rho^{-{1\over3}},\ \ \ \mu\in[1,2].
\end{equation}
Obviously,\ $\rho \gg 1 \Leftrightarrow \delta \ll 1$. Under this change of variable, the Poincar\'{e} map $P$ of (\ref{g5})  is transformed into
\begin{equation*}
\widetilde{{P}} :\left\{\begin{array}{ll}
t_1=t_0+\delta\mu_0+\delta\tilde{f}_1(\mu_0,t_0,\delta), \\[0.2cm]
\mu_1=\mu_0+\delta\tilde{f}_2(\mu_0,t_0,\delta),
\end{array}\right.
\end{equation*}
where $\tilde{f}_1,\tilde{f}_2$ are real analytic almost periodic in $t_0$ with the frequency $\omega=(\cdots,\omega_\lambda,\cdots)$ and $\tilde{f}_1,\tilde{f}_2\rightarrow 0$ as $\delta\rightarrow 0$.

Since (\ref{g13}) is a homeomorphism for $\rho$ large enough, hence from the intersection property of $P$, we know $\widetilde{P}$ also has the  intersection property. Now we apply Theorem \ref{thm6.1}  to prove  the existence of almost periodic solutions and the  boundedness of all solutions  for (\ref{g1}).

\subsection{The main results}
\begin{theorem}\label{thm7.1}
Every solution of (\ref{g1})  with a real analytic almost periodic function $f(t)\in AP_r(\omega)$, $\omega$ satisfying the nonresonance condition (\ref{b6}) is bounded. Moreover (\ref{g1}) has infinitely many almost periodic solutions.
\end{theorem}
\noindent\textbf{Proof}: From Theorem \ref{thm6.1}, we know that the Poincar\'{e} map $\widetilde{P}$ has invariant curves  if $\delta$ is sufficiently small. If the conditions of Theorem \ref{thm7.1} hold and $\delta>0$ is sufficiently small,\ therefore the assumptions of the Theorem \ref{thm6.1} are met. Hence the existence of the invariant curves of $\widetilde{P}$ is guaranteed by Theorem \ref{thm6.1},\ the invariant curves are  real analytic almost periodic  with frequencies $\omega=(\cdots,\omega_\lambda,\cdots)$.\ Undoing the change of variables we obtain the invariant curves of $P$.\ Then system (\ref{g1}) has infinitely many almost periodic solutions  as well as the boundedness of solutions.\qed

\begin{remark}\label{rem7.2}
It follows  from the proof of Theorem \ref{thm7.1} that if the conditions of Theorem \ref{thm7.1} hold,\ then system (\ref{g2}) has infinitely many almost periodic solutions with frequencies $\Big\{\omega=(\cdots,\omega_\lambda,\cdots), {1\over{\delta\alpha}}\Big\}$ satisfying the following nonresonance condition
\begin{align*}
&\Big|{\langle k,\omega \rangle {{ \delta\alpha} \over {2\pi}}-j}\Big|\geq {\gamma \over {\Delta([[k]])\Delta(|k|)}},\ \ \ \  \mbox{for all}\ \  k \in \mathbb{Z}_{0}^{{\mathbb{Z}}}\backslash\{0\},\ \  j \in \mathbb{Z},\\[0.2cm]
&\alpha\in \big[1+\gamma, 2-\gamma\big],
\end{align*}
and $\gamma, \delta$ are sufficiently small.
\end{remark}

\bigskip

\section*{References}
\bibliographystyle{elsarticle-num}

\end{document}